\documentclass{article}

\usepackage{amssymb,amsmath,theorem,euscript}

\input{epsf}

\newcounter{sec}

\newcounter{punct}[sec]

\def\punct{\refstepcounter{punct}{\arabic{sec}.\arabic{punct}.  }\boldmath}

\def\COUNTERS{\addtocounter{sec}{1}
              \setcounter{punct}{0}
          \setcounter{equation}{0}
          \setcounter{theorem}{0}
                  }

\newtheorem{theorem}{Theorem}[sec]
\newtheorem{proposition}[theorem]{Proposition}

\newtheorem{lemma}[theorem]{Lemma}

\newtheorem{corollary}[theorem]{Corollary}

\newtheorem{conjecture}[theorem]{Conjecture}

 \def\ov{\overline}
\def\wt{\widetilde}
\def\wh{\widehat}
 \newcommand{\rk}{\mathop {\mathrm {rk}}\nolimits}

 \newcommand{\indef}{\mathop {\mathrm {indef}}\nolimits}
 \newcommand{\dom}{\mathop {\mathrm {dom}}\nolimits}
  \newcommand{\End}{\mathop {\mathrm {End}}\nolimits}
   \newcommand{\Mor}{\mathop {\mathrm {Mor}}\nolimits}
    \newcommand{\Ob}{\mathop {\mathrm {Ob}}\nolimits}
     \newcommand{\Aut}{\mathop {\mathrm {Aut}}\nolimits}
 
\begin{document}

\def\OO{\mathrm{O}}
\def\GLO{\mathrm{GLO}}
\def\Coll{\mathrm{Coll}}
\def\kappa{\varkappa}
\def\Mat{\mathrm{Mat}}
\def\U{\mathrm U}
\def\GLL{\overline{\mathrm {GL}}}
\def\Spp{\overline{\mathrm {Sp}}}
\def\SL{\mathrm{SL}}
\def\PGL{\mathrm{PGL}}

\def\R{\mathbb{R}}
\def\C{\mathbb{C}}
\def\V{\mathbb{V}}

\def\la{\langle}
\def\ra{\rangle}

 \def\cA{\mathcal A}
\def\cB{\mathcal B}
\def\cC{\mathcal C}
\def\cD{\mathcal D}
\def\cE{\mathcal E}
\def\cF{\mathcal F}
\def\cG{\mathcal G}
\def\cH{\mathcal H}
\def\cJ{\mathcal J}
\def\cI{\mathcal I}
\def\cK{\mathcal K}
 \def\cL{\mathcal L}
\def\cM{\mathcal M}
\def\cN{\mathcal N}
 \def\cO{\mathcal O}
\def\cP{\mathcal P}
\def\cQ{\mathcal Q}
\def\cR{\mathcal R}
\def\cS{\mathcal S}
\def\cT{\mathcal T}
\def\cU{\mathcal U}
\def\cV{\mathcal V}
 \def\cW{\mathcal W}
\def\cX{\mathcal X}
 \def\cY{\mathcal Y}
 \def\cZ{\mathcal Z}
 
 \def\cGL{\mathcal{GL}}
 
 \def\EXP{\mathcal{EXP}}
 
 \def\frA{\mathfrak A}
 \def\frB{\mathfrak B}
\def\frC{\mathfrak C}
\def\frD{\mathfrak D}
\def\frE{\mathfrak E}
\def\frF{\mathfrak F}
\def\frG{\mathfrak G}
\def\frH{\mathfrak H}
\def\frI{\mathfrak I}
 \def\frJ{\mathfrak J}
 \def\frK{\mathfrak K}
 \def\frL{\mathfrak L}
\def\frM{\mathfrak M}
 \def\frN{\mathfrak N} \def\frO{\mathfrak O} \def\frP{\mathfrak P} \def\frQ{\mathfrak Q} \def\frR{\mathfrak R}
 \def\frS{\mathfrak S} \def\frT{\mathfrak T} \def\frU{\mathfrak U} \def\frV{\mathfrak V} \def\frW{\mathfrak W}
 \def\frX{\mathfrak X} \def\frY{\mathfrak Y} \def\frZ{\mathfrak Z} \def\fra{\mathfrak a} \def\frb{\mathfrak b}
 \def\frc{\mathfrak c} \def\frd{\mathfrak d} \def\fre{\mathfrak e} \def\frf{\mathfrak f} \def\frg{\mathfrak g}
 \def\frh{\mathfrak h} \def\fri{\mathfrak i} \def\frj{\mathfrak j} \def\frk{\mathfrak k} \def\frl{\mathfrak l}
 \def\frm{\mathfrak m} \def\frn{\mathfrak n} \def\fro{\mathfrak o} \def\frp{\mathfrak p} \def\frq{\mathfrak q}
 \def\frr{\mathfrak r} \def\frs{\mathfrak s} \def\frt{\mathfrak t} \def\fru{\mathfrak u} \def\frv{\mathfrak v}
 \def\frw{\mathfrak w} \def\frx{\mathfrak x} \def\fry{\mathfrak y} \def\frz{\mathfrak z} \def\frsp{\mathfrak{sp}}
 \def\bfa{\mathbf a} \def\bfb{\mathbf b} \def\bfc{\mathbf c} \def\bfd{\mathbf d} \def\bfe{\mathbf e} \def\bff{\mathbf f}
 \def\bfg{\mathbf g} \def\bfh{\mathbf h} \def\bfi{\mathbf i} \def\bfj{\mathbf j} \def\bfk{\mathbf k} \def\bfl{\mathbf l}
 \def\bfm{\mathbf m} \def\bfn{\mathbf n} \def\bfo{\mathbf o} \def\bfp{\mathbf p} \def\bfq{\mathbf q} \def\bfr{\mathbf r}
 \def\bfs{\mathbf s} \def\bft{\mathbf t} \def\bfu{\mathbf u} \def\bfv{\mathbf v} \def\bfw{\mathbf w} \def\bfx{\mathbf x}
 \def\bfy{\mathbf y} \def\bfz{\mathbf z} \def\bfA{\mathbf A} \def\bfB{\mathbf B} \def\bfC{\mathbf C} \def\bfD{\mathbf D}
 \def\bfE{\mathbf E} \def\bfF{\mathbf F} \def\bfG{\mathbf G} \def\bfH{\mathbf H} \def\bfI{\mathbf I} \def\bfJ{\mathbf J}
 \def\bfK{\mathbf K} \def\bfL{\mathbf L} \def\bfM{\mathbf M} \def\bfN{\mathbf N} \def\bfO{\mathbf O} \def\bfP{\mathbf P}
 \def\bfQ{\mathbf Q} \def\bfR{\mathbf R} \def\bfS{\mathbf S} \def\bfT{\mathbf T} \def\bfU{\mathbf U} \def\bfV{\mathbf V}
 \def\bfW{\mathbf W} \def\bfX{\mathbf X} \def\bfY{\mathbf Y} \def\bfZ{\mathbf Z} \def\bfw{\mathbf w}
\def\B{\mathrm{B}}
 \def\R {{\mathbb R }} \def\C {{\mathbb C }} \def\Z{{\mathbb Z}} \def\H{{\mathbb H}}
  \def\K{{\mathbb K}}
   \def\k{{\Bbbk}}
 \def\N{{\mathbb N}} \def\Q{{\mathbb Q}} \def\A{{\mathbb A}} \def\T{\mathbb T} 
 \def\G{\mathbb G}
 \def\bbA{\mathbb A} \def\bbB{\mathbb B} \def\bbD{\mathbb D} \def\bbE{\mathbb E} \def\bbF{\mathbb F} \def\bbG{\mathbb G}
 \def\bbI{\mathbb I} \def\bbJ{\mathbb J} \def\bbL{\mathbb L} \def\bbM{\mathbb M} \def\bbN{\mathbb N} \def\bbO{\mathbb O}
 \def\bbP{\mathbb P} \def\bbQ{\mathbb Q} \def\bbS{\mathbb S} \def\bbT{\mathbb T} \def\bbU{\mathbb U} \def\bbV{\mathbb V}
 \def\bbW{\mathbb W} \def\bbX{\mathbb X} \def\bbY{\mathbb Y} \def\kappa{\varkappa} \def\epsilon{\varepsilon}
 \def\phi{\varphi} \def\le{\leqslant} \def\ge{\geqslant}

\def\P{\mathbf P}

\def\GL{\mathrm {GL}}
\def\bGL{\mathbf {GL}}
\def\GLB{\mathrm {GLB}}

\def\bGr{\mathbf {Gr}}
\def\Gr{\mathrm {Gr}}
\def\Sp{\mathrm {Sp}}
\def\bFl{\mathbf {Fl}}

\def\1{\mathbf {1}}
\def\0{\mathbf {0}}

\def\rra{\rightrightarrows}

 \newcommand{\Dim}{\mathop {\mathrm {Dim}}\nolimits}
  \newcommand{\codim}{\mathop {\mathrm {codim}}\nolimits}
   \newcommand{\im}{\mathop {\mathrm {im}}\nolimits}
\newcommand{\ind}{\mathop {\mathrm {ind}}\nolimits}
\newcommand{\graph}{\mathop {\mathrm {graph}}\nolimits}

\def\F{\bbF}

\def\lambdA{{\boldsymbol{\lambda}}}
\def\alphA{{\boldsymbol{\alpha}}}
\def\betA{{\boldsymbol{\beta}}}
\def\gammA{{\boldsymbol{\gamma}}}
\def\deltA{{\boldsymbol{\delta}}}
\def\mU{{\boldsymbol{\mu}}}
\def\nU{{\boldsymbol{\nu}}}
\def\epsiloN{{\boldsymbol{\varepsilon}}}
\def\phI{{\boldsymbol{\phi}}}
\def\psI{{\boldsymbol{\psi}}}
\def\kappA{{\boldsymbol{\kappa}}}

\def\sm{\smallskip}
\def\nw{\nwarrow}
\def\se{\searrow}

\begin{center}
\Large\bf

Groups $\GL(\infty)$ over  finite fields and multiplications of double cosets

\medskip

\large \sc
Yury A. Neretin%
\footnote{Supported by the grants of FWF (Austrian science fund), P25142,P28421, P31591.}$^,$%
\footnote{The present work includes preprint 
\cite{Ner-preprint}, {\tt arxiv.org/abs/1310.1596}.}
\end{center}

{\small Let $\F$ be a finite field. Consider a direct sum $V$ of an infinite number of copies
of $\F$, consider the dual space $V^\diamond$, i.~e., the direct product of an infinite number
of copies of $\F$. Consider the direct sum $\V=V\oplus V^\diamond$.  The object of the paper is
the group $\GLL$ of  continuous linear operators in $\V$.
 We reduce
the theory of unitary representations of $\GLL$ to projective representations of
a certain category whose morphisms are  linear relations in finite-dimensional linear spaces 
over $\F$. In fact we consider a certain  family $\ov Q_\alphA$
of subgroups in $\GLL$ preserving two-element flags, show that there is a natural
multiplication on spaces of double cosets with respect to $\ov Q_\alphA$,
and reduce this multiplication to products of linear relations. 
We show that this group has type $\mathrm{I}$ and obtain an 'upper estimate' of the set of all irreducible unitary representations of  
$\GLL$.
}

\section{Formulation of results}

\COUNTERS
This section contains Introduction (Subsect. \ref{ss:introduction}),
Notation (Subsect. \ref{ss:notation}), formulation of results
(Subsect. \ref{ss:i-multiplication}--\ref{ss:i-generality}),
additional comments,  links with other works, and discussion of similar objects
(Subsect. \ref{ss:produce}--\ref{ss:poly}).

\sm

{\bf\punct Introduction.\label{ss:introduction}}
Denote by $\F$ a finite field with $q=p^l$ elements, where $p$ is prime.
Consider the linear space $\V$ over $\F$ consisting of two-sided sequences
$$
(\dots, v_{-1}, v_0,v_1,\dots),
$$
where $v_{j}=0$ for sufficiently large $j$. The space $\V$ is locally compact with respect to the natural topology%
\footnote{An infinite locally compact separable linear space over $\F$ can be discrete (such space $V$ is the direct sum  of a countable 
number of copies of $\F$), compact (such space $V^\diamond$ is isomorphic to the direct product of a countable 
number of copies of $\F$), or  isomorphic to $\V$. Groups $\GL(V)$ and $\GL(V^\diamond)$ are discussed in Subsect. \ref{ss:review}.}. So it has a unique up to a scalar factor Haar measure (i.e., a measure invariant with respect to all translations).
 Denote by $\GLL^\circ$ the group of all continuous linear operators in $\V$, by $\GLL$ its subgroup consisting of transformations preserving the Haar measure.
 Clearly $\GLL$ is a normal subgroup in $\GLL^\circ$, and 
 $$\GLL^\circ/\,\GLL\simeq \Z.$$ 
The group $\GLL$ was introduced in \cite{Ner-finite} as the maximal group of symmetries of a certain infinite-dimensional Grassmannian over a finite field.

\sm 

Now there exist well-developed theories of unitary representations of
infinite-dimensional real classical groups and of infinite symmetric groups,
these two theories are parallel one to another.
Infinite-dimensional classical groups over finite fields were topics of many attacks since \cite{Skud}, 1976,
see, e.g., \cite{Olsh-semi}, \cite{Dud}, \cite{VK1}, \cite{VK2}, \cite{GO}, \cite{GKV}, \cite{Ner-finite}, \cite{Tsa} (we present a brief survey of various versions of such groups in Subsect. \ref{ss:review}),
but a picture remains to be less transparent, less connected, and less rich than
for real groups and symmetric groups.

\sm 

 The present paper contains an attempt to classify all unitary representations of  $\GLL$. For an
 irreducible unitary representation $\rho$ of  $\GLL$ we assign a
 number $z$ satisfying $0\le z\le 1$. Next,
 
 \sm 
 
 ---
 if $z>0$, then we assign 
  a canonically defined pair $(k,\sigma)$, where $k$ is a non-negative integer and $\sigma$ is an irreducible representation of $\GL(k,\F)$;
  
  \sm 
  
 --- the case $z=0$  is slightly different, we assign two integers $l\ge m$ and an
  irreducible representation of the group $\GL(l-m,\F)$.

  \sm 
  
   Data of such type uniquely determine a representations  $\rho$, but we do not know, which $(k,\sigma)$ actually correspond to representations of $\GLL$  for $z>0$.
   
   \sm
 
 The proof is based on the following phenomenon, which is interesting by itself. Denote by $W_k\subset \V$ the subspace consisting of all vectors $v$ satisfying $v_j=0$ for $j>k$. For a pair of subspaces $W_k\supset W_l$ we denote by $\ov Q_{l,k}\subset \GLL$ the subgroup
 consisting of all $g$ such that
 
 \sm 
 
 1) $g W_k=W_k$, $g W_l= W_l$.
 
 \sm 
 
 2) $g$ induces the unit operator in the quotient space $W_k/W_l$. 

\sm 

In other words, we consider the group of all invertible  matrices of the block form
$$
\begin{pmatrix}
*&*&*\\0&1&*\\0&0&*
\end{pmatrix}.
$$

We show that {\it there is a natural associative multiplication of double coset spaces}
$$
\ov  Q_{l_1,k_1}\setminus \GLL/\,\ov  Q_{l_2,k_2}
\,\,\times \,\,
\ov Q_{l_2,k_2}\setminus \GLL/\,\ov  Q_{l_3,k_3}
\,\,\to\,\,
\ov  Q_{l_1,k_1}\setminus \GLL/\,\ov  Q_{l_3,k_3}.
$$
In this way we get a category $\cGL$, whose objects are pairs $(l,k)$ and  sets of morphisms
$(l_2,k_2)\to (l_1,k_1)$  are double coset spaces $Q_{l_1,k_1}\setminus \GLL/\ov  Q_{l_2,k_2}$.  Unitary representations of the group $\GLL$ are in a canonical one-to-one correspondence with $*$-representations of the category $\cGL$.

Next, we obtain a transparent description of the category $\cGL$.
For any double coset $\ov Q_{l_1,k_1}\cdot g\cdot\ov  Q_{l_2,k_2}$ we assign a linear subspace in $W_{k_2}/W_{l_2}\oplus W_{k_1}/W_{l_1}$ (a {\it linear relation}) and a non-negative integer. These data uniquely determine a double coset, and the category $\cGL$ is equivalent to a certain 'central extension' of the category of linear relations. This relatively easily implies the statement about representations mentioned above. 

\sm 

{\sc Remark.} Multiplication of double cosets is a fairly common phenomenon  for infinite-dimensional groups, we briefly discuss this in Subsect. \ref{ss:mult}, but our case is unusual inside a collection
of known examples.
\hfill $\boxtimes$ 

\sm

{\bf\punct Notation.%
	\label{ss:notation}}
Below $\Z_+$  denotes the set of non-negative integers, $\C^\times$ denotes the multiplicative
group of complex numbers $\C$.

By $1_k$ we denote the unit matrix of order $k$, sometimes we omit a subscript $k$, also
sometimes we write $\1$. The symbol 1 denotes also units in groups and semigroups.
 For a matrix $A$ denote by $A^t$ 
the transposed matrix. In particular, for a vector-row $v$ we denote by $v^t$
the corresponding vector-column.

Denote by $\F$ a finite field with $q=p^l$ elements, where $p$ is prime.
Denote  by $\F^\times$ its multiplicative group.

 {\it We are mainly interested in groups of infinite matrices over the finite field
$\F$}. However,
some nontrivial considerations (see Theorems \ref{th:well}, \ref{th:completeness}, \ref{th:isomorphism}) are valid in a wider generality%
\footnote{In particular, a wider generality can be interesting
in a context of $p$-adic classical groups, see \cite{Ner-p-hua}.}. Denote by
 $\fro$  an arbitrary commutative ring with unit, by $\k$  an arbitrary field.

 Denote by $\GL(\infty,\fro)$  the group of infinite matrices 
 $g=\begin{pmatrix} g_{11}& g_{12}&\dots\\
 g_{21}& g_{22}&\dots\\
 \vdots&\vdots&\ddots
 \end{pmatrix}  
 $
  over a ring $\fro$ 
  such that $g-1$ has finite number of non-zero matrix elements. We call such matrices
 $g$ {\it finitary}.

Denote by $\GL(2\infty,\fro)$ the group of all finitary two-sided infinite matrices 
$g$ over
$\fro$.
Of course, a bijection between $\N$ and $\Z$ induces an isomorphism between 
$\GL(\infty,\fro)$ and $\GL(2\infty,\fro)$.
By $\fro^{2\infty}$ we denote the space of two-sided sequences $\{v_j\}_{j\in\Z}$ consisting 
of elements of $\fro$ such that $v_j=0$ for sufficiently
large $|j|$. The group $\GL(2\infty,\fro)$ acts on this space
by multiplication of columns by matrices.

For a finite matrix $g\in\GL(\alpha,\fro)$ we define  matrices
$$
g_\se:=
\begin{pmatrix}
g&0\\
0&1_\infty
\end{pmatrix}\in \GL(\infty,\fro),
\qquad
^\nw g_\se:=
\begin{pmatrix}
1_\infty&0&0\\
0&g&0\\
0&0&1_\infty
\end{pmatrix}\in \GL(2\infty,\fro).
$$

By $S_\infty$ we denote the group of permutations of $\N$ with finite support.
By $\ov S_\infty$ we denote the group of all permutations, equipped with a natural topology
(a sequence $\sigma_j\in \ov S_\infty$ converges to $\sigma$ if for each $n\in \N$
we have $\sigma_j n=\sigma n$ starting some $j$).

\sm


Let $G$ be a group, $K$, $L$ its subgroups. By $G/\!/K$ we denote the space of conjugacy classes of $G$
with respect to $K$, by $K\setminus G/L$ the double coset space with respect to $K$ and $L$.

\sm 

A {\it Hilbert space} means a separable Hilbert space.

\sm


{\bf\punct Multiplication of double cosets.%
\label{ss:i-multiplication}}
Denote by $\cA$ the set of all pairs 
 $\alpha_-\le\alpha_+$ of integers, we denote them
by 
$$\alphA:=(\alpha_-,\alpha_+).$$
Denote
$$
|\alphA|:=\alpha_+-\alpha_-.
$$
Define a partial order on $\cA$ assuming
that 
$$
\betA\prec\alphA \qquad \text{if the segment $[\beta_-+1, \beta_+]$
is contained in $[\alpha_-+1,\alpha_+]$}. 
$$


Fix $\alphA\in\cA$.
Split $\Z$ into  segments,
\begin{equation}
-\infty<k\le \alpha_-+1,\qquad \alpha_-+1<k\le \alpha_+,\quad \alpha_+<k<\infty
\label{eq:split}
\end{equation}
of lengths $\infty$, $\alpha_+-\alpha_-$, $\infty$. 
According (\ref{eq:split}), we split our space $\fro^{2\infty}$ into a direct sum of 3 subspaces
consisting of sequences of the form
\begin{align*}
\fro_{\alphA}^-:&\,\,
 (\dots, v_{\alpha_--2}, v_{\alpha_--1}, v_{\alpha_-},\,0,0,\dots\dots\dots\dots\dots\dots\dots\dots\dots\dots.\dots.);\\
\fro_\alphA:&\,\,(\dots\dots\dots\dots\dots\dots,0,0, v_{\alpha_-+1},  v_{\alpha_-+2},\dots,  v_{\alpha_+},0,0,\dots\dots\dots\dots);\\
\fro_\alphA^+ :&\,\,(\dots\dots\dots\dots\dots\dots\dots\dots\dots\dots\dots\dots\dots,\,0,\,0,\,  v_{\alpha_-+1},  v_{\alpha_-+2},\dots).
\end{align*}

Denote by
$ Q_{\alpha_-,\,\alpha_+}= Q_\alphA$ the subgroup in $\GL(2\infty,\fro)$ consisting of block matrices
of the following form:
\begin{equation}
\left(
\begin{array}{ccc}
*&*&*\\
0&1_{|\alphA|}&*\\
0&0&*
\end{array}\right),
\label{eq:Q}
\end{equation}
the blocks correspond to the decompositions
$$
\fro^{2\infty}=\fro_{\alphA}^-\oplus \fro_{\alphA}\oplus \fro _{\alphA}^+. 
$$

We wish to show that for each $\alphA$, $\betA$, $\gammA\in\cA$
there is a natural multiplication on double coset spaces
\begin{equation*}
 Q_{\alphA}\setminus \GL(2\infty,\fro)\,/\, Q_{\betA}\,\,\,
\times\,\,\,
Q_{\betA}\setminus \GL(2\infty,\fro)\,/\, Q_{\gammA}\,\,\,
\to\,\,\,
 Q_{\alphA}\setminus \GL(2\infty,\fro)\,/\, Q_{\gammA}.
\end{equation*}
defined in the following way.
For a double coset $\fra \in Q_{\alphA}\setminus \GL(2\infty,\fro)\,/\, Q_{\betA}$
we  choose a representative
\begin{equation}
{A=\vphantom{\begin{pmatrix}
		a&b&c\\
		c&d&e\\
		e&f&g
		\end{pmatrix}}}^\nw\!\!
\begin{pmatrix}
a_{11}&a_{12}&a_{13}\\
a_{21}&a_{22}&a_{23}\\
a_{31}&a_{32}&a_{33}
\end{pmatrix}_\se
,
\label{eq:matrix}
\end{equation}
a size of the matrix has the form
$$
(N_-+|\alphA|+N_+)\times (M_-+|\betA|+M_+),
$$
where 
$$
N_--\alpha_-=M_--\beta_-,\qquad N_++\alpha_+=M_++\beta_+
$$

Notice that for any $\mu$, $\nu\ge 0$ the expression 
\begin{equation}
{\wt A=
	\vphantom{\begin{pmatrix}1_k&0&0&0&0\\
		0&a&b&c&0\\
		0&c&d&e&0\\
		0&e&f&g&0\\
		0&0&0&0&1_l
		\end{pmatrix}}}^\nw
\left(\begin{array}{cc|c|cc}
1_\mu&0&0&0&0\\
0&a_{11}&a_{12}&a_{13}&0\\
\hline 
0&a_{21}&a_{22}&a_{23}&0\\
\hline
0&a_{31}&a_{32}&a_{33}&0\\
0&0&0&0&1_\nu
\end{array}\right)_\se
\label{eq:add}
\end{equation}
determines the same element of $\GL(2\infty,\fro)$.

Consider two double cosets
$$
\fra\in
Q_{\alphA}\setminus \GL(2\infty,\fro)\,/\, Q_{\betA},\quad
\frb\in 
Q_{\betA}\setminus \GL(2\infty,\fro)\,/\, Q_{\gammA}.
$$
Choose their representatives
\begin{equation}
{A=\vphantom{\begin{pmatrix}
		a&b&c\\
		c&d&e\\
		e&f&g
		\end{pmatrix}}}^\nw\!\!
\begin{pmatrix}
a_{11}&a_{12}&a_{13}\\
a_{21}&a_{22}&a_{23}\\
a_{31}&a_{32}&a_{33}
\end{pmatrix}_\se,\qquad
{B=\vphantom{\begin{pmatrix}
		b_{11}&b_{12}&b_{13}\\
	b_{21}&b_{22}&b_{23}\\
	b_{31}&b_{32}&b_{33}
		\end{pmatrix}}}^\nw\!\!
\begin{pmatrix}
	b_{11}&b_{12}&b_{13}\\
b_{21}&b_{22}&b_{23}\\
b_{31}&b_{32}&b_{33}
\end{pmatrix}_\se
.
\label{eq:A-B}
\end{equation}
We can assume that sizes of matrices $A$, $B$
are 
$$
\Bigl(N_-+ |\alphA|+N_+\Bigr)\times \Bigl(M_-+ |\betA|+M_+\Bigr)
$$ 
and
$$
\Bigl(M_-+ |\betA|+M_+\Bigr)\times \Bigl(K_-+ |\gammA|+K_+\Bigr),
$$
otherwise we apply the transformation (\ref{eq:add}).
We wish to define a product of double cosets
$$
\fra \star \frb \in  Q_{\alphA}\setminus \GL(2\infty,\fro)\,/\, Q_{\gammA}.  
$$
For  matrices $A$, $B$ we denote
\begin{equation}
A^\circ:=
\left(\begin{array}{cc|c|cc}
1_{M_-}&0&0&0&0\\
0&a_{11}&a_{12}&a_{13}&0\\
\hline
0&a_{21}&a_{22}&a_{23}&0\\
\hline
0&a_{31}&a_{32}&a_{33}&0\\
0&0&0&0&1_{M_+}
\end{array}\right)
\label{eq:A-circ}
\end{equation}
and
\begin{equation}
B^\lozenge:=
\left( \begin{array}{cc|c|cc}
b_{11}&0&b_{12}&0&b_{13}\\
0&1_{M_-}&0&0&0\\
\hline
b_{21}&0&b_{22}&0&b_{23}\\
\hline
0&0&0&1_{M_+}&0\\
b_{31}&0&b_{32}&0&b_{33}
\end{array}
\right)
\label{eq:B-lozenge}
\end{equation}
  Then $\fra\star \frb$  
  is
 the double coset containing
the matrix
\begin{equation}
{A\star B:=}^\nw
[A^\circ B^\lozenge]_\se.
\label{eq:star}
\end{equation}


\begin{theorem}
	\label{th:well}
	{\rm a)} A double coset  $\fra\star \frb$ defined above does not depend
	on the choice of representatives $A$ and $B$.
	
	\sm
	
	{\rm b)} The $\star$-multiplication is associative, i. e.,
	for 
	\begin{multline*}
	\fra\in
	Q_{\alphA}\setminus \GL(2\infty,\fro)\,/\, Q_{\betA},\quad 
	\frb\in 
	Q_{\betA}\setminus \GL(2\infty,\fro)\,/\, Q_{\gammA},
	\\
	\frc\in 
	Q_{\gammA} \setminus \GL(2\infty,\fro)\,/\, Q_{\deltA}
	\end{multline*}
	we have
	$$
	(\fra\star\frb)\star\frc=	\fra\star(\frb\star\frc).
	$$
\end{theorem}

The proof of the theorem occupy Section \ref{s:2}.

\sm

Thus we get a {\it category} $\cGL=\cGL(\fro)$, the set $\Ob(\cGL)$ of its objects is $\cA$,
the sets of morphisms
are
$$\Mor(\betA,\alphA):=Q_\alphA\setminus \GL(2\infty,\fro)/Q_\betA.$$
We also denote by
$\End(\alphA):=\Mor(\alphA,\alphA)$ the set of endomorphisms of $\alphA$
and by $\Aut(\alphA)$ the group of automorphisms. The unit
$\1_\alphA\in \End(\alphA)$ is the double coset containing the unit matrix.
It is easy to see that  $\Aut(\alphA)\simeq \GL(|\alphA|,\fro)$.

\sm 

 The map $A\mapsto A^{-1}$ determines 
the maps $\fra \mapsto\fra^*$ of double coset spaces
$$
Q_{\alphA}\setminus \GL(2\infty,\fro)\,/\, Q_{\betA}\,\,\,\to\,\,\,
Q_{\betA}	\setminus \GL(2\infty,\fro)\,/\, Q_{\alphA}.
$$

\begin{proposition}
	\label{pr:involution}
	The maps $\fra \mapsto\fra^*$ define an {\rm involution} in the category $\cGL$, i.~e.,
	$$
	(\fra\star \frb)^*=\frb^*\star \fra^*.
	$$
\end{proposition}

The proof is contained in Subsect. \ref{ss:rewriting}.

\sm 

{\bf \punct Description of the category $\cGL(\k)$.%
\label{ss:i-description}} Now we consider
only groups of matrices over a field $\k$. In notation of the previous subsection
we replace $\fro$ by $\k$, for instance, we write $\k_\alphA$ instead of $\fro_\alphA$.
In this case we are able to present a transparent description
of the category $\cGL=\cGL(\k)$ in the language of linear relations.

Recall some simple notions.
Let $V$, $W$ be linear spaces over  $\k$.
A {\it linear relation} $P:V\rra W$ is a linear subspace in $V\oplus W$.
Let $P:V\rra W$, $Q:W\rra Y$ be linear relations. Their {\it product}
$QP:V\rra Y$ is the set of all $(v,y)\in V\oplus Y$, for which there exists
$w\in W$ such that $(v,w)\in P$, $(w,y)\in Q$.

\sm  

For a linear relation $P:V\rra W$ we define:

\sm  

--- the {\it kernel} $\ker P=P\cap V$;

\sm  

--- the {\it image} $\im P$ is the projection of $P$ to $W$ along $V$;

\sm  

--- the {\it domain of definiteness} $\dom P$ is the projection of $P$ to $V$ along $W$;

\sm  

--- the {\it indefiniteness} $\indef P=P\cap W$;

\sm  

--- the {\it rank }
\begin{multline*}
\rk P=\dim P-\dim \ker P-\dim \indef P=\\=
\dim\dom P-\dim \ker P=\dim\im P-\dim \indef P.
\end{multline*}

We also define the {\it pseudoinverse} linear relation $P^\square: W\rra V$ as the set of all pairs $(w,v)$
such that $(v,w)\in P$.

\sm 

{\sc Remark.} A graph of a linear operator $T:V\to W$ is a linear relation,
in this case $\dom A=V$, $\indef A=0$. A product of operators is a special case
of products of linear relations. If an operator is invertible, then
its pseudoinverse linear relation is the graph of the inverse
operator. \hfill $\boxtimes$

\sm

For an element $A$ given by (\ref{eq:matrix}) of $\GL(2\infty,\k)$ we define
a {\it characteristic linear relation}
$$
\chi(A): \k_{\betA}\rra \k_{\alphA}
$$
as the set of all pairs $(v,u)\in
\k_{\betA}\oplus \k_{\alphA}$,
for which there exist $x\in \k_{\alphA}^-$, $y\in \k_{\betA}^-$ such that
\begin{equation}
\begin{pmatrix}
x\\u\\0
\end{pmatrix}
=
\begin{pmatrix}
a_{11}&a_{12}&a_{13}\\
a_{21}&a_{22}&a_{23}\\
a_{31}&a_{32}&a_{33}
\end{pmatrix}
\begin{pmatrix}
y\\v\\0
\end{pmatrix}.
\label{eq:def-chi}
\end{equation}

\begin{lemma}
	\label{l:chi}
	The linear relation $\chi(A)$ depends only on the double coset	
	$\fra\in Q_{\alphA}\setminus \GL(2\infty,\k)\,/\, Q_{\betA}$
	containing $A$.
\end{lemma}

A double coset is not determined by the characteristic relation
and we need
 an additional invariant of a double coset, namely
$$
\eta(\fra):=\rk a_{31}.
$$

 

\begin{proposition}
	\label{pr:equality}
	Any pair  $(\chi(\fra),\eta(\fra))$ satisfies the condition
	\begin{equation}
	\eta(\fra)\ge \beta_--\alpha_-+\dim \ker \chi(\fra)-\dim \indef \chi(\fra).
	\label{eq:chi-eta}
	\end{equation}
	This is a unique restriction for such pairs.
\end{proposition}

\begin{theorem}
	\label{th:completeness}
	The pair $\bigl(\chi(\fra), \eta(\fra)\bigr)$ completely determines a double coset $\fra$.	
\end{theorem}

\begin{theorem}
\label{th:isomorphism}
	Let $\fra\in Q_\alphA\setminus\GL(2\infty,\k)\,/\,Q_\betA$,
	$\frb\in Q_\betA\setminus\GL(2\infty,\k)\,/\,Q_\gammA$. Then
	
	\sm 
	
	{\rm a)} $\chi(\fra\star \frb) =\chi(\fra)\chi(\frb)$.
	
	\sm 
	
	
	
	{\rm b)}  $\eta(\fra\star\frb)=\eta(\fra)+\eta(\frb)+
	\dim \bigl( \indef \chi(\frb)/(\indef \chi(\frb)\cap \dom \chi(\fra))\bigr)$.
	
	\sm
	
	{\rm c)} We have $\chi(\fra ^*)=\chi(\fra)^\square$ and
	\begin{equation}
	\eta(\fra^*)=\eta(\fra)+ \dim\indef \chi-\dim\ker\chi-\beta_-+\alpha_.
	\label{eq:fra*}
	\end{equation}
\end{theorem}

{\sc Remark.} The category of double cosets admits a more natural description in  terms of polyhomomorphisms, see below Subsect. \ref{ss:poly}. 
 \hfill $\boxtimes$

\sm 

Proofs of statements of this subsection are contained in Section \ref{s:3}.

\sm

{\bf\punct The group $\GLL$.%
\label{ss:i-GLL}} Now we again reduce our generality  and  
consider groups of matrices over
a finite field $\F$. Our next purpose is to define the topological group $\GLL$,
which is the main object of the paper.

Denote by $V$ the  direct sum of a countable
number of copies 
of the field $\F$.
We regard $V$ as the space of all vectors
$$
v=(v_1, v_2,v_3,\dots),
$$
where $v_j\in \F$ and $v_j=0$ for all but a finite number of $j$. 
The group $V$ is  countable, we equip it with the discrete topology.
By $V^\diamond$ we denote the direct product of an infinite number of copies of $\F$.
We regard $V^\diamond$ as the space of all vectors
$$v^\diamond=(\dots, v_{-2}, v_{-1}, v_0),$$
where $v_j\in \F$.
We equip $V^\diamond$ with the product topology
and get  a compact group. The groups $V$ and $V^\diamond$ are Pontryagin dual
(on the Pontryagin duality, see, e.g., \cite{Mor} or \cite{Kir}, Subsect. 12.1).
More precisely, 
there is a pairing $V\times V^\diamond\to \F$,
defined by
\begin{equation}
S(v,w)=\sum_{j=0}^\infty  v_{j+1} w_{-j}.
\label{eq:Svw}
\end{equation}
(a sum actually is finite).
Choosing  an arbitrary nontrivial homomorphism $\chi$ of the additive group of $\F$ to $\C^\times$
 we write the
 duality $V\times V^\diamond\to \C^\times$
by
$$(v,w)\mapsto \chi\bigl\{  S(v,w) \bigr\}.$$

The direct sum 
$$
\V:=V^\diamond\oplus V
$$
can be regarded as the space of all two-sided sequences
\begin{equation}
(\dots, v_{-2}, v_{-1},v_0,v_1,v_2,\dots),
\label{eq:two-side}
\end{equation}
where $v_j\in \F$ and $v_{m}=0$ for sufficiently large $m$.
The group $\V$ 
is a locally compact Abelian group.

Denote by $\GLL^\circ$ the group of all continuous linear transformations%
\footnote{If the order of $\F$ is prime, then $\GLL^\circ$ is the group of all automorphisms
of the Abelian group $\V$.} of $\V$,
for details, see \cite{Ner-finite}. We write elements of $V^\diamond\oplus V$ as columns,
matrices act by multiplications from the left:
\begin{equation}
\begin{pmatrix}v^\diamond\\ v\end{pmatrix}\mapsto
g \begin{pmatrix}v^\diamond\\ v\end{pmatrix}=  \begin{pmatrix}
a&b\\c&d
\end{pmatrix}
\begin{pmatrix}v^\diamond\\ v\end{pmatrix}.
\label{eq:abcd}
\end{equation}
In this notation:

\sm

1) the matrix $b:V\to V^\diamond$ can be arbitrary;

\sm

2) $c:V^\diamond\to V$  contains only a finite number of nonzero elements;

\sm

3) each row of $a:V^\diamond\to V^\diamond$ contains only a finite number nonzero elements; 

\sm

4) each column
of $d:V\to V$ contains a finite number of nonzero elements.

\sm

Also $g$ has an inverse matrix satisfying the same properties.

\sm

{\sc Remark.} We can reformulate the list 1)--4)  in the following
way. We decompose $W$ as a sum of two subspaces $\V=W_k\oplus  W^\circ$,
  $W_k$ consists of vectors $v$ such that $v_j=0$ for $j\ge k$,  $W^\circ_k$  consists 
of vectors $v$ such that $v_j=0$ for all $j\le k$.  
Represent $g$ as a block matrix according the decomposition $\V=W_k\oplus  W^\circ_k$,
$$
g=\begin{pmatrix}
A_k&B_k\\ C_k&D_k
\end{pmatrix}.
$$
Then for any $k$ the block $C_k$ contains only a finite number of nonzero matrix elements.
\hfill $\boxtimes$

\sm 

We equip the group $\GLL^\circ$ with the {\it topology} of uniform convergence on compact sets%
\footnote{Topologies on such groups are determined by  groups themselves.
	More precisely, the statement {\it Any homomorphism of Polish groups is continuous} is compatible with the Zermelo--Fraenkel system plus the axiom of dependent choice. 
	Any reasonable complete topology invented by the reader for $\GLL(\V)$  will be the same.} 
(since $V^\diamond\oplus V$ is an Abelian topological group, the uniform convergence 
is well defined, see e.g., \cite{Kir}, Subsect. 2.6, or \cite{Mor}, Sect. 8). Denote by $W_m$ the subgroup in $\V$
consisting of all sequences (\ref{eq:two-side}) such that $v_{j}=0$ for $j>m$.
We get an exhausting sequence 
$$
\dots \subset W_k\subset W_{k+1}\subset W_{k+2}\subset \dots
$$  
of compact subsets (subgroups)
in $\V$. A sequence $g_j\in \GLL^\circ$ converges to 1, if for any $\alphA=(\alpha_-,\alpha_+)$
for sufficiently large $j$ for each $w\in W_{\alpha_+}$ we have $(g_j-1)w\in W_{\alpha_-+1}$.

Equivalently, we consider the subgroups $\ov Q_\alphA$ consisting of block matrices
of the form 
$\left(
\begin{array}{ccc}
*&*&*\\
0&1_{|\alphA|}&*\\
0&0&*
\end{array}\right)$
as above, see (\ref{eq:Q}).
Then such subgroups form a fundamental system of open neighborhoods of unit. The closure of
$Q_\alphA$ is $\ov Q_\alphA$. 

\sm 

Normalize a Haar measure on $\V$ by the assumption: the measure of
$V^\diamond$ is 1.
We prefer to work with a subgroup 
$$\GLL\subset\GLL^\circ$$ 
consisting of transformations preserving the Haar measure on $\V$.
Generally, an element  $g\in \GLL^\circ$ sends a set of measure $\alpha$ 
to a set of measure $p^k\alpha$, where $k=k(g)\in \Z$, see a general statement
\cite{Bour},  Subsect. VII.1.4.

Below we need the following a more precise description of the subgroup $\GLL$
(for details, see \cite{Ner-finite}).
For a matrix $a=:a(g)$ in (\ref{eq:abcd})  there is the following analog of the {\it Fredholm index}
(see \cite{Ner-faa}, Subsects. 2.4--2.7)
defined by
$$
\ind a(g):= \codim \im a(g)-\dim\ker a(g)
$$
(both numbers are finite).
It is easy to verify that $g\mapsto \ind a(g)$ is a homomorphism $\GLL^\circ\to\Z$,
$$
\ind a(g_1g_2)= \ind a(g_1)+\ind a(g_2).
$$
We define the group $\GLL$ as the kernel of this homomorphism.

Consider the shift  operator 
\begin{equation}
J:\{v_j\}\mapsto \{v_{j+1}\}
\label{eq:shift}
\end{equation}
 in the space of sequences
(\ref{eq:two-side}).
It
is  contained
in $\GLL^\circ$ but not in $\GLL$, in fact we have a semidirect product
$$
\GLL^\circ\simeq \Z\ltimes \GLL,
$$
where $\Z$ is the group  generated by the shift.


\sm 

{\bf\punct Multiplicativity.%
\label{ss:i-multiplicativity}} 
Here we show that unitary representations of the group $\GLL$ can be reduced 
to representations of the category $\cGL(\F)$. 

\begin{lemma}
	\label{l:repr-finite}
{\rm a)}	Any double coset $\fra\in \ov Q_\alphA \setminus \GLL/\ov Q_\betA$ 
	contains an element of $\GL(2\infty,\F)$.
	
	\sm 
	
{\rm b)} If $A_1$, $A_2\in \GL(2\infty)$ are contained
in one double coset $\fra \in \ov Q_\alphA \setminus \GLL/\ov Q_\betA$,
then they are contained in one double coset	$\in  Q_\alphA \setminus \GL(2\infty,\F)/\ Q_\betA$.
\end{lemma}

Therefore the following sets coincide:
$$
\ov Q_\alphA \setminus \GLL/\ov Q_\betA\simeq
 Q_\alphA \setminus \GL(2\infty,\F)/\ Q_\betA.
$$

Let $\rho$ be a unitary representation of $\GLL$ in a Hilbert space $H$.
Denote by $H_\alphA\subset H$ the subspace of
$\ov Q_{\alphA}$-fixed vectors.
Obviously,
$$
\Bigl\{\betA\succ\alphA\Bigr\}
\Rightarrow \Bigl\{H_{\betA}\supset H_{\alphA}\Bigr\}.
$$
Denote by $P_\alphA$ the operator of orthogonal projection
$H\to H_\alphA$.

\begin{lemma}
	\label{l:admissibility}
	For a unitary representation $\rho$ of $\GL(2\infty,\F)$ in a Hilbert
	space $H$ the following statements are equivalent:
	
	\sm 
	
	{\rm 1)} $\rho$ has a continuous extension to the group $\GLL$;
	
	\sm 
	
	{\rm 2)} the subspace $\cup_{\alphA} H_{\alphA}$ is dense
	in $H$.
\end{lemma}

For $\fra\in \ov Q_\alphA \setminus \GLL/\ov Q_\betA$ we 
define the operator
$$
\wh\rho_{\alphA,\betA}(\fra): H_\betA\to H_\alphA
$$
by
$$
\wh\rho_{\alphA,\betA}(\fra):= P_\alphA \rho(A)\bigr|_{H_\beta},
\qquad
\text{where $A\in \fra$.}
$$
It can be readily checked that this operator does not depend
on the choice of a representative $A\in\fra$.
The following {\it multiplicativity theorem} holds:

\begin{theorem}
	\label{th:multiplicativity}
	{\rm a)} For any $\alphA$, $\betA$, $\gammA\in\cA$ and any
	$$
	\fra\in \ov Q_\alphA \setminus \GLL/\ov Q_\betA,
	\qquad \frb \in \ov Q_\betA \setminus \GLL/\ov Q_\gammA
	$$
	we have
	$$
	\wh\rho_{\alphA,\betA}(\fra)\, \wh\rho_{\betA,\gammA}(\frb)=
	\wh\rho_{\alphA,\gammA}(\fra\star \frb).
	$$
	
	{\rm b)} 	For any $\alphA$, $\betA$ and any
	$
	\fra\in \ov Q_\alphA \setminus \GLL/\ov Q_\betA$
	we have
	\begin{equation}
	\wh\rho_{\alphA,\betA}(\fra)^*=	\wh\rho_{\betA,\alphA}(\fra^*).	
	\label{eq:*}
	\end{equation}
	
	{\rm c)} $\|\wh\rho_{\alphA,\betA}(\fra)\|\le 1$.
\end{theorem}

The statement a) requires a proof, b) and c) are obvious.

\sm 

In other words we get a representation of the category $\cGL(\F)$,
i.e., a functor from the category $\cGL(\F)$ to the category
of Hilbert spaces and bounded linear operators (the  first statement
of the theorem). This representation
is compatible with the involution
(we also use the term {\it $*$-representations}
for  such representations., i.e., representations satisfying
(\ref{eq:*})).

\sm 

Proofs of statements of this subsection are contained in Section \ref{s:4}.

\sm

{\bf \punct The spherical character of irreducible representations of $\GLL$.%
\label{ss:i-spherical-character}}
Here for any unitary representation of $\GLL$ in a Hilbert space $H$ we assign a canonical
 intertwining  operator $\frZ:H\to H$.
 By the Schur lemma, for an irreducible representation $\rho$
this operator is scalar, $\frZ=z\cdot 1$. We call the number $z=z(\rho)$ by the {\it spherical character
of $\rho$.}

Denote by $\zeta_\alphA^k$ the double coset defined by
\begin{equation}
{\zeta_\alphA^k=
	\ov Q_\alphA\cdot
	\vphantom{\begin{pmatrix}
		1&d\\0&\sigma''
		\end{pmatrix}}}^\nw 
\begin{pmatrix}
0&0&1_k\\
0&1_{|\alphA|}&0\\
1_k&0&0
\end{pmatrix}_\se
\cdot \ov Q_\alphA. 
\label{eq:zeta}
\end{equation}

The following statement is straightforward:

\begin{proposition}
	{\rm a)}	
	The element $\zeta_\alphA^k$ is contained in the center
	of the semigroup $\End(\alphA)=\ov Q_\alphA\setminus \GLL/\,\ov Q_\alphA$ and
	$$
	\zeta_\alphA^k\star \zeta_\alphA^l=\zeta_\alphA^{k+l}.
	$$
	
	{\rm b)} For each $\alphA$, $\betA$ and $\fra\in \Mor(\betA,\alphA)$ we have
	$$
	\zeta_\alphA^k \star \fra= \fra\star \zeta_\betA^k.
	$$
	
	{\rm b)} The characteristic linear relation $\chi(\zeta_\alphA^k)$ 
	is the graph of the unit operator, $\eta(\zeta_\alphA^k)=k$.
\end{proposition}

The statement b) means that for $\alphA\prec\betA$ we have
$$
\wh \rho_{\betA,\betA}(\zeta_\betA^k)\Bigr|_{H_\alphA}=\wh \rho_{\alphA,\alphA}(\zeta_\alphA^k)
.$$ 
 Therefore we have  well defined self-adjoint operators 
 $\frZ^k$ in $H$ satisfying
 $$
 \frZ^k \Bigr|_{H_\alphA}=\wh \rho_{\alphA,\alphA}(\zeta_\alphA^k),\qquad
 \|\frZ^k\|\le 1,\qquad \frZ^k \frZ^l=\frZ^{k+l}
 $$ 
 commuting with the representation $\rho$. In particular,
 for  irreducible
 representations such operators $\frZ^k$ must be scalar operators,
 $$
 \frZ^k h=z^kh,\qquad \text{where $-1\le z\le 1$.}
 $$
 
 \begin{lemma}
 	\label{l:z}
 The number $z(\rho)$ satisfies to the condition	$0\le z\le 1$.
 \end{lemma}

The proof of the lemma is given in Subsect. \ref{ss:z}.

\begin{conjecture}
	\label{con:1}
	The number $z(\rho)$ has the form $p^{-k}$, where $k\in \Z_+$, or $z=0$.
	\end{conjecture}

{\sc Remarks.} 
a) There is an explicit compact semigroup $\Gamma\subset G$  with separately continuous
multiplication, such that $G$ is dense in $\Gamma$ and
any unitary representation admits a weakly continuous extension to $\Gamma$, see 
  \cite{Ner-polihomo}. This semigroup semigroup has a center isomorphic
  to  $\Z_+$, the center acts 
  in unitary representations of $G$ by operators $\frZ^k$.
  
  \sm 
  
 b) An existence of such 'spherical character' is a general phenomenon for infinite-dimensional
 groups, see \cite{Ner-faa}, Subsect. 1.12, \cite{Ner-symm}, Subsect. 3.17.
 However, usually these characters depend on an infinite number of parameters. 
\hfill $\boxtimes$.

\sm 

{\bf \punct Data determining irreducible representations of the group $\GLL$.%
	\label{ss:i-generality}}
Thus Theorem \ref{th:multiplicativity} reduces unitary representations 
of the group $\GLL$ to representations of the category $\cGL(\F)$.
Here we discuss some corollaries from this reduction.

\begin{proposition}
\label{th:trivial}
Let $\rho$ be an irreducible unitary representation
of $\GLL$ in a Hilbert space $H$.
 If $H_\alphA\ne 0$, then $\rho$ is uniquely
 determined by the  representation
 of $\End(\alphA)$ in $H_\alpha$ {\rm(}which is automatically 
irreducible{\rm )}.
 \end{proposition}

This is an automatic statement in the spirit of: 'an irreducible unitary representation $\sigma$
of a group $G$ is determined by any its matrix element $\la \sigma(g) v,v\ra$'.
See \cite{Ner-symm}, e.g., Lemma 2.7.

 For a representation $\rho$ denote by
 $\Xi(\rho)$ the set of all $\alphA$ such that $H_\alphA\ne0$.
 Clearly, if $\alphA\in \Xi(\rho)$ and $\betA\succ \alphA$,
 then $\betA\in\Xi(\rho)$.
 
 \begin{lemma}
 	\label{l:zero-nonivertible}
Let $\alphA$ be a minimal element of $\Xi(\rho)$. Let for
$\fra\in \End (\alphA)$ we have $\wh \rho(\fra)\ne 0$.
Then%
\footnote{Formally, we must say $\chi(\fra)$ is a graph
	of an invertible operator.}
 $\chi(\fra)\in \GL(|\alphA|,\F)$.
 \end{lemma}

If additionally $\rho$ is irreducible, then we have an irreducible  representation of
 of the semigroup
$\GL(|\alphA|)\times \Z_+$.  It is determined by
an irreducible representation of $\GL(|\alphA|)$ and a number $z\in[0,1]$, defining
a representation  of $\Z_+$.

 \begin{theorem}
 		\label{th:induced}
 	 Let $\rho$ be an irreducible unitary representation of $\GLL$.
 
 \sm 
 
 {\rm a)} Let $z(\rho)>0$. Then there is $k\ge 0$ such that set $\Xi(\rho)$ consists of all
 $\alphA$ such that $|\alphA|\ge k$. If $|\alphA|=|\alphA'|=k$,
 then the representations of $\GL(k,\F)$ in $H_\alphA$ and $H_{\alphA'}$ are equivalent.
 
 \sm 
 
 {\rm b)} Let $z(\rho)=0$. Then there is $\alphA$ such that $\Xi(\rho)$ consists 
 of all $\betA\ge \alphA$.
\end{theorem} 

\begin{figure}
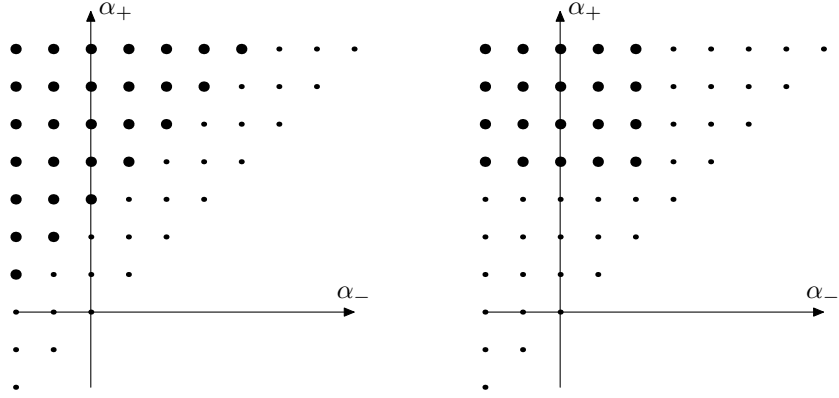

	
		$$\epsfbox{support.1}\qquad\qquad \epsfbox{support.2}$$
	
	\caption{To Theorem \ref{th:induced}
		 A set $\Xi(\rho)$ for $z>0$ and $z=0$.}
	 \label{fig:Xi}
		
\end{figure}

See, Fig. \ref{fig:Xi}.

Thus, for $z>0$ an irreducible representation $\rho$ of $\GLL$ is uniquely determined
by a triple $(z,k,\tau)$, where $\tau$ is an irreducible representation
of $\GL(k,\F)$. For $z=0$ a representation $\rho$ is determined by
a triple $(0,\alphA,\tau)$, where $\tau$ is
 an irreducible representation  of $\GL(|\alphA|,\F)$.

\sm

{\sc Remark.} Let $J:\V\to\V$ be the shift operator, see \ref{eq:shift}.
Then the map $g\mapsto JgJ^{-1} $ is an automorphism of the group $\GLL$.
In particular for any irreducible representation $\rho(g)$ of $\GLL$
we have a representation $\rho^J(g):=\rho(JgJ^{-1})$. If $z\ne0$, then
$\rho^J$ is equivalent to $\rho$. It is easy to show  that $\rho$ can be extended
to a unitary representation of the group $\GLL^\circ$.
If $z=0$, then $\rho^J$ is not equivalent to $\rho$.
The set $\Xi(\rho^J)$ 
 (see Fig. \ref{fig:Xi}.b) is obtained from  $\Xi(\rho)$
by the shift $(\alpha_-,\alpha_+)\mapsto (\alpha_-+1,\alpha_++1)$
\hfill $\boxtimes$

\begin{corollary} \label{pr:sphericity}
	The subgroup $\ov Q_{0,0}$ is spherical%
	\footnote{A subgroup $H$ in a topological group $G$ is {\it spherical} 
		if for any irreducible unitary representation of $G$ the dimension of the subspace
		of $H$-fixed vectors is $\le 1$. The definition assumes that $H$-fixed vectors exist 
		in some representations of $G$. In our case this is so, for instance the 
		natural representation of $\GLL$ in $\ell^2(\V)$ has a $Q_{0,0}$-fixed vector,
		see also \cite{Ner-finite}.} in $\GLL$.
\end{corollary}

Indeed, the corresponding semigroup $\End(\mathbf{0})=\GL(0,\F)\times \Z_+$
is $\Z_+$, its irreducible representations are one-dimensional and a spherical representation
is determined by its spherical character.
\hfill $\square$

 \sm

Return to a general situation. The case $z=0$ is rather simple. 
Consider the subgroup $\ov P_\alphA\supset\ov Q_\alphA$ consisting of matrices
of the form
$$
\begin{pmatrix}
a_{11}&a_{12}&a_{13}\\
0& a_{22}&a_{23}\\
0&0&a_{33}
\end{pmatrix}.
$$ 
 Then $\ov Q_\alphA$ is normal in $\ov P_\alphA$,
 and
 $$
 \ov P_\alphA/\,\ov Q_\alphA\simeq \GL(|\alphA|,\F).
 $$
 A representation $\tau$ of $\GL(\alphA,\F)$ can be regarded
 as a representation of $\ov P_\alphA$. It is easy to see
 that the homogeneous space $\GLL/\,\ov P_\alphA$ is countable,
 therefore the induction of representations from
 $P_\alphA$ to $\GLL$ makes sense.
 
 \begin{proposition}
 	\label{pr:z=0}
 The representation of $\GLL$ corresponding to a triple
 $(0,\alphA,\tau)$ is the representation induced from
 the representation $\tau$ of the subgroup $\ov P_\alphA$. 	
 \end{proposition}
 
Representation theory of finite groups $\GL(n,\F)$ (see, e.g., \cite{Macd}, Chapter IV)
does not appear in our considerations. It seems that it must  be reflected in a final
form of a classification of representations%
\footnote{and in problems of harmonic analysis discussed in the next subsection.}.

\begin{conjecture}
For a given $k$ and an irreducible representation
$\tau$ of $\GL(k,\F)$ the set of $z>0$ such that the 
 triple $(z,k,\tau)$ corresponds to a unitary representation
 of $\GLL$ has the form $p^{-l}$, where $l$ is integer
 ranging in a set of form $l\ge m$, where $m=m(k,\tau)$.	
	\end{conjecture}

\begin{theorem}
\label{th:I}
{\rm a)}The group $\GLL$ has%
\footnote{For the definition and discussion of types of representation and groups, see, e.g., 
	\cite{Mac}, \cite{Dix}.}
 type $I$.
 
 \sm 
 
{ \rm b)}
  Any irreducible unitary representation
of $\GLL$ is a direct integral  of irreducible representations. 
\end{theorem}

{\sc Remark.} If Conjecture \ref{con:1} is the truth,
then any unitary representation of $\GLL$ is a direct {\it sum} 
of  irreducible representations (this is clear from the proof in Subsect. \ref{ss:proof-I}).
\hfill $\boxtimes$

\sm

Statements of this subsection are proved in Section \ref{s:5}.

At this point we end the presentation  of theorems of the paper,
our next purpose is  additional remarks
and some links, formal and heuristic, with other works on infinite-dimensional groups.

\sm 

{\bf\punct 
	Some problems of harmonic analysis for the group $\GLL$.%
\label{ss:produce}}
In this paper we do not discuss constructions of irreducible representations
of the group $\GLL$, except  Subsect. \ref{ss:constructions}.
In any case, there are the following problems of harmonic analysis, which provide us a way
to obtain a lot of irreducible representations in their spectra.

\sm

{\sc 1) Howe duality.} Let the characteristic of $\F$ be $\ne 2$.
 Define a skew-symmetric bilinear form on $\V=V^\diamond\oplus V$
by
$$
\bigl\{  (v_1,w_1), (v_2,w_2)\bigr\}:=S(v_1,w_2)-S(v_2,w_1),
$$
where  $S(\cdot,\cdot)$ is the pairing $V^\diamond\times V\to\F$, see (\ref{eq:Svw}).
Denote by $\Spp=\Spp(\V)$ the subgroup of $\GLL$ consisting of operators in $\V$ preserving this form.
For this group the Weil representation is well defined: since $V$ and $V^\diamond$ are locally
compact,  the construction of Weil
\cite{Wei} remains to be valid in this case. As it was shown in \cite{Ner-Weil},
this representation admits an extension to a certain category of 'Lagrangian'
linear relations $\V\rra\V$.

The group $\GLL$ admits a natural embedding to $\Spp$. Namely,
 consider the space $
 \V:=V^\diamond\oplus V$ and its dual
$$\V^\diamond=(V^\diamond\oplus V)^\diamond\simeq V\oplus V^\diamond\simeq \V.$$
For $A\in \GLL$  consider the dual operator $A^t$ in $\V^\diamond$ and its inverse $(A^t)^{-1}$.
The group $\GLL$ acts in $\V\oplus \V^\diamond$ by operators
$$
A\mapsto \begin{pmatrix}
          A&0\\0&(A^t)^{-1}
         \end{pmatrix}
$$
preserving the duality $\V\times \V^\diamond\to \F$.
On the other hand
$$\V\oplus \V^\diamond
=(V^\diamond\oplus V)\oplus (V\oplus V^\diamond)
\simeq (V^\diamond\oplus V^\diamond)\oplus (V\oplus V)\simeq\V,$$
since $V\oplus V\simeq V$.
We get the desired embedding $\GLL(\V)\to \Spp(\V\oplus \V)$.

Next, fix $m$, consider the space $(\V\oplus \V)\otimes \F^m$,
 consider the Weil representation of $\Spp\bigl((\V\oplus \V)\otimes \F^m\bigr)$
and restrict it to the subgroup  $\GL(\V)\times \GL(m,\F)$.
We come to a question of the Howe duality type (see, e.g., \cite{Howe2}, \cite{GW}). 
The Howe duality for groups over finite fields was a topic of numerous works
(see, e.g., \cite{AMR}). Notice that  finite-dimensional counterparts of our objects
are not  pairs $(\GL(N,\F), \GL(m,\F))$ but $(\End_\cGL(\alphA), \GL(m,\F))$.
 May be our problem  is more similar to the initial Howe's work
\cite{Howe}.

\sm 

2) {\sc Flag Spaces.}
 In \cite{Ner-finite} there were constructed  $\GLL$-invariant
measures on certain
spaces of flags in $\V$. Namely, there were considered (finite or infinite) flags of discrete cocompact%
\footnote{I.e., quotients $\V/L_{i}$ are compact.}
subspaces 
$$\dots \subset L_{j-1}\subset L_j\subset L_{j+1}\subset \dots $$
in $\V$.  
There arises a problem about explicit decompositions of the corresponding spaces $L^2$.
For the Grassmannians this problem has an explicit solution in terms of Carlitz--Al Salam 
$q$-hypergeometric orthogonal polynomials, the spectrum consists of $Q_{0,0}$-spherical representations.

\sm

3) {\sc Other spaces of flags.} We can also consider finite flags
$$
M_1\subset \dots \subset M_k
$$ 
of compact subspaces in $\V$ of positive Haar measure. Such spaces of flags  are countable
and Mackey's argumentation \cite{Mack-induced} (see also \cite{Cor})  reduces
the decomposition of $\ell^2$ to certain questions about finite groups
$\GL(n,\F)$.
It seems that this
 problem  is less interesting than 1) and 2).

\sm

{\bf\punct General remarks on multiplication of conjugacy classes and double cosets.\label{ss:mult}} Here we discuss some standard
facts related to classical groups.
Let $\k$ be a field and $\fro\subset\k$ be a subring with unit
(the most important case is%
\footnote{The case when $\k$ is a $p$-adic field
	$\Q_p$ and $\fro$ is a ring of $p$-adic integers was discussed in
	\cite{Ner-rational}.}
$\k=\fro=\C$). 
Denote by $\GL(m+\infty,\k)$
the same group $\GL(\infty,\k)$ considered as a group of finitary block matrices
$\begin{pmatrix}a&b\\c&d \end{pmatrix}$
of size $m+\infty$. Denote by $K$ the subgroup, consisting of matrices
of the form
$\begin{pmatrix}1_m&0\\0&H \end{pmatrix}$, where $H$ is a matrix over $\fro$,
$$K\simeq \GL(\infty,\fro).$$
We claim that the set of conjugacy classes
$$G/\!/K=\GL(m+\infty,\k)/\!/\GL(\infty,\fro)$$
is a semigroup with respect to the following $\circ$-multiplication.
For two matrices 
$$\begin{pmatrix}a&b\\c&d \end{pmatrix}_\se,
\qquad \begin{pmatrix}p&q\\r&t \end{pmatrix}_\se\in \GL(m+\infty)$$
we define their $\circ$-product by
\begin{multline}
\begin{pmatrix}a&b\\c&d\end{pmatrix}_\se\circ
\begin{pmatrix}p&q\\r&t\end{pmatrix}_\se
:=
\left[\begin{pmatrix}a&b&0\\c&d&0\\0&0&1 \end{pmatrix}
\begin{pmatrix}p&0&q\\0&1&0\\r&0&t\\ \end{pmatrix}\right]_\se
=\\=
\left(
\begin{array}{c|cc}ap&b&aq
\\ 
\hline
cp&d&cq\\r&0&t \end{array}\right)_\se
.
\label{eq:circ}
\end{multline}

\begin{theorem}
	\label{th:colligations}
	{\rm a)}
	The $\circ$-multiplication is a well-defined associative operation on the set
	of conjugacy classes
	$\GL(m+\infty,\k)/\!/\GL(\infty,\fro)$.
	
	\sm
	
	{\rm b)} The $\circ$-multiplication is a well-defined associative operation on the set
	of conjugacy classes
	$\U(m+\infty)/\!/\U(\infty)$. 
\end{theorem}

The statement (as soon as it is formulated) is more-or-less obvious.
Various versions of these  semigroups are classical topics of system theory and operator theory,
see, e.g., \cite{Bro}, \cite{Dym}, Chapter 19, \cite{GGK}, Part VII, 
\cite{Haz}, \cite{RR}, see also \cite{Ner-rational}. 

\sm

If $\fro=\k$, then
the multiplication $\circ$ can be clarified in the following way.
For $g=\begin{pmatrix}a&b\\c&d \end{pmatrix}$ we
write the following%
\footnote{Cf. the definition (\ref{eq:def-chi}) of the characteristic linear relation.}
 'perverse equation for eigenvalues':
\begin{equation}
\begin{pmatrix}
p\\x
\end{pmatrix}=
\begin{pmatrix}a&b\\c&d \end{pmatrix}
\begin{pmatrix}
q\\ \lambda x
\end{pmatrix}, \qquad \lambda\in \k.
\end{equation}
Eliminating $x$ we get
a relation of the type
$$p=\chi_g(\lambda) q,$$
where $\chi: \k\mapsto \GL(m,\k)$
is the '{\it Livshits characteristic function}' or '{\it transfer-function}'
$$
\chi_g(\lambda)=a+\lambda b(1-\lambda d)^{-1}c.
$$

\begin{theorem}
	$$\chi_{g\circ h} (\lambda)= \chi_{g} (\lambda)\, \chi_{h} (\lambda).$$
\end{theorem}

\sm

{\sc Multiplications of double cosets.}
 We denote:

\sm

--- $\U(\infty)\subset \GL(\infty,\C)$ is
the group of finitary unitary matrices over $\C$; 

\sm

--- $\OO(\infty)\subset \GL(\infty,\R)$ is  the group
of finitary real orthogonal  matrices;

\sm 

---  $\OO(\infty,\C)\subset \GL(\infty,\C)$ is
the group of finitary complex orthogonal matrices.

\begin{theorem}
	\label{th:colligations-2}
	{\rm a)}
	The $\circ$-multiplication is a well-defined associative operations on double coset
	spaces
	$\GL(\infty,\fro)\setminus\GL(m+\infty,\k)/\,\GL(\infty,\fro)$.
	
	\sm
	
	{\rm b)} The formula {\rm(\ref{eq:circ})} determines an associative operation on double cosets
	\begin{equation}
	\OO(\infty)\setminus \U(m+\infty)/\,\OO(\infty).
	\label{eq:OUO}
	\end{equation}
	
	
	{\rm c)} The formula {\rm(\ref{eq:circ})} determines an associative operation on double cosets
	$\OO(\infty,\C)\setminus \GL(m+\infty,\C)/\,\OO(\infty,\C)$.
\end{theorem}

In the case of $\OO(\infty)\setminus \U(m+\infty)/\,\OO(\infty)$
we have  a 'multiplicativity theorem' as we discussed above. 
See \cite{Olsh-GB}, for details, see \cite{Ner-book},  Section IX.4.



 According
\cite{Olsh-GB}, theorems of
this type hold for all infinite-dimensional limits of symmetric pairs $G\supset K$.
Recently \cite{Ness1}, \cite{Ness2}, \cite{Ner-char}, \cite{Ner-faa}, \cite{Ner-char}, \cite{Ner-symm}
it was observed that these phenomena are quite general. For instance, 
there is a well defined multiplication on the  double cosets space 
\begin{equation}
\mathrm{diag}\, \U(\infty)\setminus
\underbrace{
 \GL(m+\infty,\C)\times \dots\times\GL(m+\infty,\C)}_{\text{$m$ times}}/\,\mathrm{diag}\, \U(\infty),
\label{eq:times}
\end{equation}
where $\mathrm{diag}\, \U(\infty)$ is the subgroup 
$$
\mathrm{diag} \,\U(\infty)\subset \mathrm{diag}\, \GL(\infty,\C)\subset\mathrm{diag}\, \GL(m+\infty,\C)
$$ 
in the diagonal $\mathrm{diag}\, \GL(m+\infty,\C)$ of the direct product.
Such multiplication can be described in terms
of  semigroups of matrix-valued rational functions
of matrix argument.

 In all known  cases we have double cosets 
 with respect to  infinite dimensional analogs $K$ of certain simple or reductive 
groups as
 \begin{equation}
 \U(n), \quad \U(n)\times \U(n), \quad \mathrm{O}(n),\quad \Sp(n),\quad S_n,\quad \dots
 \end{equation}

It seems natural to continue this list by $\GL(n,\F)$. However,  we consider 
 groups of block matrices over $\F$ having the form
\begin{equation}
\begin{pmatrix} *&*&*\\0&1&*\\0&0&*\end{pmatrix}
\label{eq:***}.
\end{equation}
This family of  subgroups naturally arose in the context of \cite{Ner-finite}.
Lemma \ref{l:two-subspaces} below provides us an a priori explanation:  any vector in a representation
of $\GLL$ fixed by all matrices of the form
$\begin{pmatrix} *&0&0\\0&1&0\\0&0&*\end{pmatrix}$
is also fixed by all matrices of the form (\ref{eq:***}).

Another unexpected  place is the proof
of Theorem \ref{th:well}.a (which claims that the product of double cosets
is well defined). The proof is not difficult, but meet an obstacle, which is not observable
in previously known cases, see below Remark after Lemma \ref{l:well}.
In particular, I do not see a possibility to produce 
counterparts of semigroups (\ref{eq:times}) considering double cosets 
with respect to subgroups of the type (\ref{eq:***}).

\sm

{\bf \punct What is an infinite-dimensional group $\GL$ over a finite field?%
\label{ss:review}}
The group
$$
\GL(\infty, \F)=\lim\limits_{\longrightarrow} \GL(n,\F)
$$
of finitary matrices
is not%
\footnote{The simplest proof: the group $\PGL(\infty,\F)=\GL(\infty,\F)/\F^\times$
has infinite conjugacy classes except the unit. This easily implies that 
the left regular representation of $\PGL(\infty,\F)$ generates a Murray-von Neumann factor
of the type $\mathrm{II}_1$, see \cite{Mac}, Corollary on page 62. On the hand
by the Thoma criterion a discrete group has type I iff  it has an Abelian subgroup 
of finite index.}
 of type I. So a representation theory in the usual sense for this group is impossible. There were several approaches to formulate  problems of representation theory for this group.

 \sm

{\sc A. Unitary representations of completions.}
Define the following completions of $\GL(\infty,\F)$:

\sm 

1) The group $\GLL(V)$ of all linear operators in the space $V$ (recall that $V$ is the direct
sum of a countable number of copies of $\F$). In other words, we consider the group of all invertible
matrices $g$ such that $g$ and $g^{-1}$ have only finite number of nonzero matrix elements
in each column.

\sm

2) The group $\GLL(V^\diamond)$ of all linear operators in the space $V^\diamond$ (recall that $V^\circ$ is the direct
product of a countable number of copies of $\F$).  The groups $\GLL(V^\diamond)$ and $\GLL(V)$ are isomorphic,
the isomorphism is given by $g\mapsto (g^t)^{-1}$.

\sm 

3) The group 
$$
\GLL(V\sqcup V^\diamond):= \GLL(V)\cap \GLL(V^\diamond).
$$
In this case $g$ and $g^{-1}$ have only finite number of nonzero matrix elements in each 
column and each row%
\footnote{
	 For completeness, we say definitions of topologies.  The group $\GLL(V)$ acts by permutations of a countable set
	$V$, its topology is induced from  the symmetric group. The group $\GLL(V\sqcup V^\diamond)$ 
	has two different actions on $V$, namely, $v\mapsto gv$ and $v\mapsto (g^t)^{-1}v$. So it acts
	on a countable space
	 $V\sqcup V$, the topology also is induced from the symmetric group.}.

\sm

 Denote by
$G_\alpha$, $Q_\alpha$, $P_\alpha$ the subgroups in $\GL(\infty,\F)$ consisting of block matrices of size $\alpha+\infty$ having the form
$$
\begin{pmatrix}
1&0\\0&*
\end{pmatrix},
\qquad\qquad
\begin{pmatrix}
1&*\\0&*
\end{pmatrix},\qquad
\qquad
\begin{pmatrix}
*&*\\0&*
\end{pmatrix} \qquad
 \text{respectively}.
$$
For a unitary representation of $\GL(\infty,\F)$ in a Hilbert space
$H$ we denote by $H^{G_\alpha}$ (resp. $H^{Q_\alpha}$) the subspaces in $H$ consisting of $G_\alpha$-fixed (resp.  $Q_\alpha$-fixed) vectors. 

\sm 

Classification of unitary representations of the group $\GL(V)$ was obtained in 2012
by Tsankov\cite{Tsa} as a special case of his general theorem on oligomorphic groups. The answer
is simple: any irreducible representation is induced from a representation of
a subgroup $P_\alpha$ (for some $\alpha=0$, 1, 2, \dots) trivial on $Q_\alpha$.

It can be shown (this is a simplified version of Theorem \ref{th:well})
that double cosets 
$$Q_\alpha \setminus \GL(\infty,\F)/Q_\beta$$
form a category and $*$-representations 
of this category  are in a one-to-one correspondence with unitary representations of $\GLL(V)$. Explicit description of this category is simple: this is the category
of partial isomorphisms%
\footnote{A {\it partial isomorphism} of linear spaces $R:X\to Y$
	is a bijection of a subspace in $X$ to a subspace in $Y$. In other words it is a linear relation with $\ker R=0$, $\indef R=0$}  of finite-dimensional linear spaces over $\F$.

\sm 

The same work of Tsankov covers the group $\GL(V\sqcup V^\diamond)$, but in a certain sense
unitary representations of this group were described by Olshanski 
 \cite{Olsh-semi}, 1991. Olshanski  considered the class of unitary representations  of
the group $\GL(\infty,\F)$  {\it admissible} in the following sense:
the subspace
$
\cup_{\alpha} H^{G_\alpha}
$
is dense in $H$. The paper \cite{Olsh-semi} contains a classification of all admissible representations (see, also Dudko \cite{Dud}).

\begin{proposition}
	\label{pr:ts}
	A unitary representation of $\GL(\infty,\F)$ is Olshanski admissible if and only if 
	it admits a continuous extension to $\GLL(V\sqcup V^\diamond)$.
\end{proposition}

Olshanski noted the arrow $\Rightarrow$. Tsankov observed the arrow $\Leftarrow$,
which follows from coincidence of  classifications of admissible representations of
$\GL(2\infty,\F)$  and representations of $\ov\GL(V\sqcup V^\diamond)$.
%
For a clarification of the picture, we present an a priori proof  in Subsect. 
\ref{ss:digression}.

\sm 

In this case 
double cosets 
$$G_\alpha\setminus \GL(\infty,\F)/G_\beta$$
form a category (this is a special case of Theorem \ref{th:colligations-2}.a), admissible representations of $\GL(\infty,\F)$ are in one-to-one correspondence with $*$-representations of this category.
This category has not faithful $*$-representations%
\footnote{The same phenomenon arises  for infinite-dimensional
$p$-adic groups, see \cite{Ner-p-adic-2}.}.
 Consider the common kernel of all  representations
and the quotient category by the kernel. This leads to the following category.
Denote by $Y_n$ the space $\F^n$, by $Y_n^\diamond$ the dual space. Denote by
$\{\cdot,\cdot\}$ the pairing $Y\times Y^\diamond\to\F$. An object of the category
is the direct sum $Y_n\oplus Y_n^\diamond$. A morphism 
 $Y_n\oplus Y_n^\diamond\to Y_m\oplus Y_m^\diamond $ is  a pair of partial isomorphisms
 $\sigma:Y_n\to Y_m$, $\sigma^\diamond:Y_n^\diamond\to Y_m^\diamond$
 such that 
 $$
 \{\sigma(z), \sigma^\diamond(z^\diamond) \}=\{z, z^\diamond \},\qquad \text{where $z\in \dom \sigma$,\,
 $z^\diamond\in \dom \sigma^\diamond$}.
 $$

Our group $\GLL$ contains a subgroup  $\GLL(V)\times \GLL(V^\diamond)$. If to look to the analogy with real classical groups
(see \cite{Ner-book}), it seems that $\GLL(V)$ is a counterpart of heavy groups and
$\GLL$ is a counterpart of the Olshanski infinite-dimensional classical groups. 
On the other hand the set of parameters of representations of $\GLL$ seems small comparatively
infinite dimensional real groups or infinite symmetric groups.

\sm 

{\sc B.  Infinite-dimensional Hecke algebras.}
Thoma  \cite{Tho-symm}, 1964, classified all representations of the infinite symmetric group $S_\infty$ generating  Murray--von Neumann factors%
\footnote{For definitions, see, e.g., \cite{Mac} or \cite{Dix}.} of type $\mathrm{II}_1$. This is equivalent to a description of extreme points of the set of central positive definite functions on $S_\infty$. Olshanski \cite{Olsh-symm} noticed that this problem also is equivalent to a description of representations of the double
$S_\infty\times S_\infty$ spherical with respect to the diagonal $S_\infty$.
Skudlarek \cite{Skud} in 1976 tried to extend the Thoma approach to $\GL(\infty,\F)$
but his list of positive definite central functions was trivial. 

\sm 

Define some groups and subgroups:

---  denote by 
  $\GLB(\infty,\F)$ the group  of all  matrices having only finite
number of nonzero
elements under the diagonal; this completion of $\GL(\infty,\F)$
is the next topic of our overview;

\sm

 --- denote by $\mathrm{B}(\infty,\F)\subset \GLB(\infty,\F)$
  the group consisting of upper triangular matrices.
 
 \sm 
 
 --- let  $\GLB(n,\F)\subset \GLB(\infty,\F)$
  be the group generated by $\GL(n,\F)$ and
 $\B(\infty,\F)$, i.e., the group of all infinite invertible matrices
 $g_{ij}$ such that $g_{ij}=0$ whenever $i>j$, $i>n$.
 
 \sm 
 
  Then $\GLB(\infty,\F)$ is the inductive limit
 $$
 \GLB(\infty,\F)=\lim_{\longrightarrow} \GLB(n,\F).
 $$
 The groups $\B(\infty,\F)$, $\GLB(n,\F)$
 are  compact, $\mathrm{GLB}(\infty,\F)$
  is locally compact and  is not a group
of type $\mathrm{I}$.

For a locally compact  group $G$ and its compact open subgroup
$K$ denote by $\cA(K\setminus G/K)$ the convolution algebra consisting
of compactly supported continuous functions, which are constant on double cosets
$K\cdot g\cdot K$. 
In other words, we consider the algebra of $K$-biinvariant functions $f$
on $G$: for $k_1$, $k_2\in K$, we have $f(k_1 g k_2)=f(g)$.

If $G=\GL(n,\F_q)$ and $K=\B(n,\F_q)$ is the group of upper triangular matrices,
then 
$$\cA\bigl(\B(n,\F_q)\setminus \GL(n,\F_q)/\B(n,\F_q)\bigr)$$
 is the well-known {\it Hecke--Iwahori
algebra} $\cH_q(n)$ of dimension $n!$, see, Iwahori, \cite{Iwa}, 1964.
It is generated by double cosets $s_j:=\B(n,\F_q)\sigma_k \B(n,\F_q)$, where
$\sigma_j\in \GL(n,\F_q)$ is the permutation of $j$-th and $(j+1)$-th
basis elements in $\F^n$, relations
are
\begin{align}
s_i s_j&= s_j s_i; \qquad \text{if $|i-j|\ge 2$;}
\label{eq:iwa1}
\\
s_j s_{j+1} s_j&=s_{j+1} s_j s_{j+1};
\label{eq:iwa2}
\\
s_j^2&=(q-1)s_j+q e,
\label{eq:iwa3}
\end{align}
 where $e$ is the double coset $\B(n,\F_q)\cdot 1\cdot \B(n,\F_q)$.
 We also have an antilinear involution defined by $\sigma_j^*=s_j$,
 $(ab)^*=b^*a^*$.
 Initially, $\cH_q$ was defined for $q=p^l$ being a power of prime.
 But relations allow to consider this algebra for any $q\in \C$.
 Clearly, for $q=1$ this algebra $\cH_1(n)$ is the group algebra of the symmetric group $S_n$.
 
  On the other hand (see \cite{GKV}, Proposition 2.5)
$$
\cA\bigl(\B(n,\F_q)\setminus \GLB(n,\F_q)/\B(n,\F_q)\bigr)\simeq
\cA\bigl(\B(n,\F_q)\setminus \GL(n,\F_q)/\B(n,\F_q)\bigr).
$$ 
We have inclusions 
$$
\B(\infty,\F)\setminus \GLB(n,\F)/\,\B(\infty,\F)\subset 
\B(\infty,\F)\setminus \GLB(n+1,\F)/\,\B(\infty,\F),
$$
this allows to regard the algebra of $\GLB(\infty,\F)$-biinvariant
functions on the group $ \GLB(\infty,\F)$ as 
the inductive limit
\begin{multline}
\cA(\mathrm{B}(\infty,\F)\setminus \GLB(\infty,\F)/\,\mathrm{B}(\infty,\F)\bigr)
= \lim_{\longrightarrow} 
\cA\bigl(\B(\infty,\F)\setminus \GLB(n,\F)/\,\B(\infty,\F)\bigr)
=\\=
\cup_{n=1}^\infty 
\cA\bigl(\B(\infty,\F)\setminus \GLB(n,\F)/\,\B(\infty,\F)\bigr).
\label{eq:limit}
\end{multline}
This algebra $\cH_q(\infty)$ is generated by $s_1$, $s_2$, \dots,
relations are given by the same formulas
(\ref{eq:iwa1})--(\ref{eq:iwa3}).

\sm 

Vershik and Kerov in 1988 \cite{VK0} obtained a  classification of all extreme positive traces
(extreme traces of $\cH_q(n)$ correspond to irreducible  characters of $\GL(n,\F_q)$)
 on $\cH_q(\infty)$, for $q>0$ (a trace $T$ is positive if $T(aa^*)\ge 0$). 
If  $q=1$, then the classification coincides with the Thoma's
classification for $S_\infty$. 


Some further works%
\footnote{S.~V.~Kerov died in 2000, Vershik published  text \cite{VK1} based on his posthumous notes.} are \cite{VK1}, \cite{GKV}; in \cite{GKV} it was  announced  that
any extreme positive trace on $\cH_q$ generates an irreducible unitary representation of the double 
$\GLB(\infty,\F_q)\times \GLB(\infty,\F_q)$; the restriction of such a representation to 
a single  $\GLB(\infty,\F_q)$ generates a Murray-von Neumann factor  of type $\mathrm{I}$ or $\mathrm{II}_\infty$.

 
 
 
 \sm 
 
 {\sc C. Few words about a comparison.}
 It seems (at least in the present moment), that stories with $\GLB(\infty,\F)$
 and $\GLL$ are orthogonal. In any case, both of them are based on limits 
 of algebras of the type 
 $$\cA\bigl (K(n)\setminus G(n)/K(n)\bigr), \qquad
 \text{where $G(1)\subset G(2)\subset\dots$,\,
 $K(1)\subset K(2)\subset\dots$.}
 $$
 Rather often, a limit algebra naturally degenerates to a semigroup
 (or to a semigroup algebra of a semigroup). A mechanism of degeneration
 is
 explained in \cite{Olsh-semi}, \cite{Ner-conv}. In a certain sense, for sufficiently large values of $n$ a convolution
 of uniform measures concentrated by  given double cosets
 $K(n)g_1 K(n)$ and  $K(n)g_2 K(n)$
 is concentrated near a third double coset  $K(n)g_3 K(n)$.
 On the other hand a double coset $K(\infty) g K(\infty)$ generates a
 well-defined operator in the space of $K(\infty)$-fixed vectors
 in a unitary representation, and quite often the set of such operators is closed
 with respect to multiplication%
  \footnote{Firstly, this phenomenon arose in the work by Ismagilov,
 	\cite{Ism}, he considered $G=\SL(2,\K)$, where $\K$ is a complete normed field
 	having infinite ring of residues. The subgroup $K$ is $\SL(2)$ over integers
 	of the field. The semigroup of double cosets in this case is $\Z_+$.}.
 The inductive limit (\ref{eq:limit}) for $\GLB(\infty,\F)$ is unusual in the following
 sense:
 only groups $G(n)=\GLB(n,\F)$ change, the prelimit compact
 subgroups $K(n)=\B(\infty,\F)$ remain to be the same.
 For this reason, algebras $\cA(\dots)$ range into an inductive limit
 (\ref{eq:limit}).
 
 A degeneration of a product simplifies a situation. On the other hand,
 this allows to enrich picture, since we can include to consideration
 objects of the type (\ref{eq:OUO}). (\ref{eq:times}) or numerous
 examples in \cite{Ner-symm}. In such cases prelimit objects seem to be
 unapproachable (at least in the present moment).

 \sm 
 
 {\sc D. On the group $\mathrm{GLB}(\infty,\F)$ and the Steinberg representation.}
The Steinberg representations of the group $\mathrm{GLB}(\infty,\F)$ and a similar group
$\mathrm{GLB}(2\infty,\F)$ were considered in \cite{GKV} and \cite{Ner-steinberg},
the second paper contains an explicit realization of the Steinberg representation
in a certain reproducing kernel Hilbert space. The Steinberg representation is a delicate
and fascinating construction of representation theory of groups $\GL(n,\F)$. Its surviving
in $\mathrm{GLB}(\infty,\F)$-limit begs the  question
about existence of  certain interesting natural class of irreducible representations
of  $\mathrm{GLB}(\infty,\F)$.

\sm

{\bf \punct Infinite-dimensional Chevalley groups.%
\label{ss:chevalley}} There are the following groups, for which our approach must work,
at least partially.




1) The symplectic group $\Spp(2\infty)$ defined in Subset. \ref{ss:produce}.

\sm 

2) The orthogonal group $\ov{\mathrm{O}}(2\infty)$ of the space $\V$, i.e.,
the group of operators in $\V$ preserving the bilinear form
$$
\bigl[  (v_1,w_1), (v_2,w_2)\bigr]:=S(v_1,w_2)+S(v_2,w_1).
$$
We also can add a one dimensional summand to this space and get an  infinite orthogonal group
of 'odd order' $\ov{\mathrm{O}}(2\infty+1)$.

\sm 

3) If $q=p^{2l}$, then the field $\F$ has an automorphism of order 2, namely
$$ 
x\mapsto \ov x:=x^{p^l}.
$$
In this case we also have the 'unitary' group $\ov \U(2\infty,\F_q)\subset \GLL$
consisting of matrices preserving the sesquilinear form
$$
\bigl[  (v_1,w_1), (v_2,w_2)\bigr]:=S(v_1,\ov w_2)+S(w_1, \ov v_2).
$$

\sm

{\bf\punct  The category $\cGL(\F)$ as a category of polyhomomorphisms.%
\label{ss:poly}} 
Here we  present another interpretation of invariants $\chi(\fra)$, $\eta(\fra)$,
the inequality (\ref{eq:chi-eta}) for these invariants, and
the formula (\ref{eq:fra*}).

Normalize  a Haar measure on $\V$ assuming that the measure of $W_0$ is 1.
Our spaces $\F_\alphA$ are quotients $W_{\alpha_+}/W_{\alpha_-}$. The Haar measure 
determines a uniform measure $\mu_\alphA$ on each quotient: the measure of each point is 
$p^{\alpha_-}$, the total measure of $W_{\alpha_+}/W_{\alpha_-}$ is $p^{\alpha_+}$.
Fix $\alphA$, $\betA$.
Let $A\in \GLL$. Consider the subspace $Y=A^{-1}W_{\alpha_+}\cap W_{\beta_-}\subset \V$
and consider the map 
$$S:Y\to W_{\beta_+}\oplus W_{\alpha_+}$$ given by
$
S:y\mapsto(y, Ay).
$
Passing to quotients, we get a map
$$
\sigma:Y\to W_{\beta_+}/W_{\beta_-}\oplus W_{\alpha_+}/W_{\alpha_+}.
$$
This a rephrasing of  equation (\ref{eq:def-chi}) determining the characteristic linear relation. So $\sigma(Y)=\chi(\fra)$.
But the space $Y$ is also equipped with a Haar measure, its image under
$\sigma$ is a canonically defined uniform measure $\nu_\fra$ on the subspace $\chi(\fra)$.
A measure of each point is $p^{\beta_- -\eta(\fra)- \dim\indef\chi(\fra)}$,
\begin{align*}
\Bigl\{\text{projection of $\nu_\fra$ to $\F_\betA$}\Bigr\}=p^{-\eta(\fra)}\mu_\betA\Bigr|_{\dom \chi(\fra)};
\\
\Bigl\{\text{projection of $\nu_\fra$ to $\F_\alphA$}\Bigr\}=p^{-\eta(\fra^*)}\mu_\alphA\Bigr|_{\im \chi(\fra)}.
\end{align*}
On this language, the passage $\fra\mapsto\fra ^*$ is simply the permutation
$\F_\betA\oplus \F_\alphA\to\F_\alphA\oplus \F_\betA$. We also see that projections of $\nu_\fra$ are
dominated by $\mu_\alphA$ and $\mu_\betA$ (and this explains the inequality (\ref{eq:chi-eta})).

\sm

Let $G_1$, $G_2$ be locally compact groups equipped with fixed two-side invariant Haar measures
$dg_1$, $dg_2$ respectively. 
According \cite{Ner-polihomo} a {\it polyhomomorphism} $(H,dh):(G_1,dg_1)\to (G_2,dg_2)$
is a closed subgroup $H\subset H_1\times H_2$ with fixed Haar measure
$dh$ such that projection of $dh$ to $G_1$ (resp. to $G_2$) is dominated by $dg_1$ (resp. by
$dg_2$). For polyhomomorphisms $(H,dh):(G_1,dg_1)\to (G_2,dg_2)$,
$(K,dk):(G_2,dg_2)\to (G_3,dg_3)$ there is a well-defined product and the 
product obtained in Theorem \ref{th:isomorphism} is a special case
of the product of polyhomomorphisms.

\sm 

{\bf \punct The further structure of the paper.} 
In Section \ref{s:2} we prove Theorem 1.1 about products of double cosets.
The description of this product in terms of linear relations is derived in
Section \ref{s:3}.  Multiplicativity is proved in Section \ref{s:4}.
The statements (Theorem \ref{th:induced}--\ref{th:I}) on representations of 
$\GLL$ are obtained in Section \ref{s:5}.

\section{Multiplication of double cosets\label{s:2}}

\COUNTERS

Here we prove  the statements of Subsect. \ref{ss:i-multiplication},
i.e.,  we show that  the category of double cosets is well defined.

\sm

{\bf\punct Rewriting  of the definition.%
\label{ss:rewriting}}
Recall, see (\ref{eq:star}), that
\begin{equation}
A\star B{=\vphantom{B}}^\nw [A^\circ B^\lozenge]_\se,
\label{eq:circ-lozenge}
\end{equation}
where $A^\circ$ is defined by (\ref{eq:A-circ}) and $B^\lozenge$ by (\ref{eq:B-lozenge}).
Denote
\begin{equation}
J_\betA(\nu,\mu)=
\left(
\begin{array}{cc|c|cc}
0&1_\nu&0&0&0\\
1_\nu&0&0&0&0\\
\hline
0&0&1_{\betA}&0&0\\
\hline
0&0&0&0&1_\mu\\
0&0&0&1_\mu&0
\end{array}\right)
\label{eq:Jm}
.
\end{equation}
We have
$$
B^\lozenge=J_{\betA}(M_+,M_-) B^\circ J_{\betA}(M_+,M_-)
$$
Therefore we can rewrite (\ref{eq:circ-lozenge}) as
$$
A\star B{=\vphantom{B}}^\nw [A^\circ J_{\betA}(M_+,M_-) B^\circ J_{\betA}(M_+,M_-)]_\se.
$$
Since $J_\betA(M_+,M_-)\in Q_\betA$, the same double coset is given by
the formula
\begin{equation} 
Q_\alphA\cdot \,
{\vphantom{B}}^\nw [A^\circ J_{\betA}(M_+,M_-) B^\circ ]_\se\,\cdot Q_\gammA
\label{eq:symmetry}
\end{equation}
This implies Proposition \ref{pr:involution} about the involution.

\sm

{\bf \punct Proof of Theorem \ref{th:well}.a.}
Denote the expression in the square brackets in (\ref{eq:symmetry})
by 
$$
A\circledast B=A^\circ J_{\betA}(M_+,M_-) B^\circ.
$$
It is sufficient to prove the following statement

\begin{lemma}
	\label{l:well}
	Let $A$, $P\in \GL(\infty,\fro)$
	 and $\Phi\in Q_\betA$. Then

	\sm
	
	{\rm a)} There exists $\Gamma\in Q_\gammA$ such that 
	\begin{equation}
	(A\cdot \Phi)\circledast P=(A\circledast P)\cdot\Gamma.
	\label{eq:factor}
	\end{equation}
	
	\sm
	
	{\rm b)} There exists $\Delta\in Q_\alphA$ such that 
	$A\circledast(\Phi\cdot P)=\Delta\cdot(A\circledast P)$.
\end{lemma}

{\sc Remark.} In all statements of this type known earlier, the factor
$\Gamma$ in (\ref{eq:factor}) depends only on $\Phi$. In our case
this factor depends on $\Phi$ and $P$, see (\ref{eq:uvw}), (\ref{eq:QVY}).
\hfill $\boxtimes$

\sm 

By the symmetry in formula (\ref{eq:symmetry}), it is
sufficient to prove the first statement.
To avoid subscripts (as in (\ref{eq:A-B})) consider matrices
\begin{equation}
A=
\begin{pmatrix}
a&b&c\\
d&e&f\\
g&h&j
\end{pmatrix}
,\quad
P=
\begin{pmatrix}
p&q&r\\
u&v&w\\
x&y&z
\end{pmatrix}
\end{equation}
 Then
$$
A\circledast P=
\left( 
\begin{array}{cc|c|cc}
p&0&q&0&r\\
bu&a&bv&c&bw\\
\hline 
eu&d&ev&f&ew\\
\hline
hu&g&hv&j&hw\\
x&0&y&0&z
\end{array}
\right).
$$

It is sufficient to prove the lemma for $\Phi$ ranging
in a collection of generators
\begin{equation}
\begin{pmatrix}
\mu&0&0\\
0&1&0\\
0&0&\nu
\end{pmatrix}
,\quad
\begin{pmatrix}
1&\phi&0\\
0&1&0\\
0&0&1
\end{pmatrix}
,\quad
\begin{pmatrix}
1&0&\theta\\
0&1&0\\
0&0&1
\end{pmatrix}
,\quad
\begin{pmatrix}
1&0&0\\
0&1&\psi\\
0&0&1
\end{pmatrix}
\label{eq:generators}
\end{equation}
of the group $Q_\betA$. Here the size of matrices is $M_-+|\betA|+M_+$, the matrices
$\mu$, $\nu$ range in $\GL(M_{\mp},\fro)$; matrices $\phi$, $\psi$, $\theta$
are arbitrary matrices (of appropriate size).

We examine these generators case by case. 

First,
$$
\left[A \cdot \begin{pmatrix}
\mu&0&0\\
0&1&0\\
0&0&\nu
\end{pmatrix} \right] \circledast P=\!\!
\begin{pmatrix}
p&0&q&0&r\\
bu&a\mu&bv&c\nu&bw\\
eu&d\mu&ev&f\nu&ew\\
hu&g\mu&hv&j\nu&hw\\
x&0&y&0&z
\end{pmatrix}
=\!\! (A\circledast P) \!\!
\begin{pmatrix}
1&0&0&0&0\\
0&\mu&0&0&0\\
0&0&1&0&0\\
0&0&0&\nu&0\\
0&0&0&0&1
\end{pmatrix},
$$
the $(K_-+M_-+|\gammA|+M_++K_+)$-matrix in the right-hand side is contained in the subgroup $Q_\betA$.

Second, 
\begin{multline}
\left[A\cdot \begin{pmatrix}
1&\phi&0\\
0&1&0\\
0&0&1
\end{pmatrix} \right]\circledast P
=
\begin{pmatrix}
p&0&q&0&r\\
bu+a\phi u &a&b+a\phi v&c&bw+a\phi w\\
eu+d\phi u&d &e+d\phi v&f&ew+d\phi w\\
hu+g\phi u&g &h+g\phi v&j&hw+g\phi w\\
x&0&y&0&z
\end{pmatrix}
=\\=
(A\circledast P) \cdot
\begin{pmatrix}
1&0&0&0&0\\
\phi u&1&\phi v&0&\phi w\\
0&0&1&0&0\\
0&0&0&1&0\\
0&0&0&0&1\\
\end{pmatrix}
\label{eq:uvw}
,\end{multline}
the $(K_-+M_-+|\gammA|+M_++K_+)$-matrix in the right-hand side is contained in the subgroup $Q_\betA$.

Next,
\begin{multline*}
\left[A\cdot \begin{pmatrix}
1&0&\theta\\
0&1&0\\
0&0&1
\end{pmatrix} \right]\circledast P
=
\begin{pmatrix}
p&0&q&0&r\\
bu&a&bv&c+a\theta &bw\\
eu&d&ev&f+d\theta &ew\\
hu&g&hv&j+g\theta &hw\\
x&0&y&0&z
\end{pmatrix}
=\\=
(A\circledast P)
\cdot
\begin{pmatrix}
1&0&0&0&0\\
0&1&0&\theta&0\\
0&0&1&0&0\\
0&0&0&1&0\\
0&0&0&0&1\\
\end{pmatrix}.
\end{multline*}

Examine the last generator of the list (\ref{eq:generators}). We have
\begin{equation}
\left[A\cdot \begin{pmatrix}
1&0&0\\
0&1&\psi\\
0&0&1
\end{pmatrix} \right]\circledast P
=
\begin{pmatrix}
p&0&q&0&r\\
bu&a&bv&c+b\psi &bw\\
eu&d&ev&f+e\psi &ew\\
hu&g&hv&j+h\psi  &hw\\
x&0&y&0&z
\end{pmatrix}.
\label{eq:psi}
\end{equation}
Denote
$$
\begin{pmatrix}
p&q&r\\
u&v&w\\
x&y&z
\end{pmatrix}^{-1}
=
\begin{pmatrix}
P&Q&R\\
U&V&W\\
X&Y&Z
\end{pmatrix}
.$$
Then the right-hand side of (\ref{eq:psi}) is
\begin{equation}
\begin{pmatrix}
p&0&q&0&r\\
bu&a&bv&c&bw\\
eu&d&ev&f&ew\\
hu&g&hv&j&hw\\
x&0&y&0&z
\end{pmatrix}\cdot
\begin{pmatrix}
1&0&0& Q\psi&0\\
0&1&0&0&0\\
0&0&1&V\psi&0\\
0&0&0&1&0\\
0&0&0&Y\psi&1\\
\end{pmatrix}.
\label{eq:QVY}
\end{equation}
(the first factor is $A\circledast P$). To verify this, we must evaluate 4th column
of the product. We get
\begin{align*}
&pQ+qV+rY=0;\\
&buQ\psi+bvV\psi+ c+bwY\psi=c+b(uQ+vV+wY)\psi=c+b\psi;\\
&euQ\psi+evV\psi+ f+ewY\psi=f+e(uQ+vV+wY)\psi=f+e\psi;\\
&huQ\psi+hvV\psi+ j+hwY\psi=j+h(uQ+vV+wY)\psi=j+h\psi;\\
&xQ+yV+zY=0,
\end{align*}
and this completes the proof.

\sm

{\bf\punct Associativity.} A group $Q_\alphA$ contains a product $S_\alphA$ of two copies of
the symmetric group $\ov S(\infty)$, it consists of 0-1-matrices
of the form
$$
\begin{pmatrix}
u&0&0\\
0&1_{|\alphA|}&0\\
0&0&v
\end{pmatrix}.
$$

To verify the statement b) of Theorem \ref{th:well}, we must show that
for  representatives $A$, $B$, $C$ of cosets $\fra$, $\frb$, $\frc$
there exist matrices $\Pi\in Q_\alphA$, $\Gamma\in Q_\delta$ such that
$$
A\star (B\star C)=\Pi\cdot (A\star B)\star C\cdot \Gamma.
$$
It is more-or-less clear that we can chose desired $\Pi\in S_\alphA$, $Q\in S_\gammA$.

\section{Description of the category of double cosets%
\label{s:3}}

\COUNTERS

Here we prove the statements of Subsect. \ref{ss:i-description}.
The proof of completeness of the system of invariants
(Theorem \ref{th:completeness}) is relatively long. We observe that the group
$\GL(|\alphA|,\F)\times \GL(|\betA|,\F)$ acts in both
the double coset space $Q_\alphA\setminus \GL(\infty,\F)/Q_\betA$
and in the target space (a linear relation plus an invariant $\in \Z$).
Next we observe that $\GL\times\GL$-orbits in two spaces are in one-to-one correspondence
and show that stabilizers of orbits coincide.

Proof of Theorem \ref{th:isomorphism} (isomorphism of categories
of double cosets and of extended linear relations) is
 parallel to proofs of previously known statements in this spirit in \cite{Ner-book}, Sect. IX.4,
and in \cite{Ner-char}.

\sm 

{\bf \punct The characteristic linear relation.%
\label{ss:invariants-proofs}} 
Here we prove that {\it the characteristic linear relation
$\chi(\cdot)$ is an invariant of a double coset $\fra$} (Lemma \ref{l:chi}).

 We consider an element $A\in \GL(2\infty,\k)$,
and write the corresponding equation (\ref{eq:def-chi}) for another element of the same double coset,
\begin{equation}
\begin{pmatrix}
x'\\u\\0
\end{pmatrix}
=
\begin{pmatrix}
d_{11}&d_{12}&d_{13}\\
0&1&d_{23}\\
0&0&d_{33}
\end{pmatrix}^{-1}
\begin{pmatrix}
a_{11}&a_{12}&a_{13}\\
a_{21}&a_{22}&a_{23}\\
a_{31}&a_{32}&a_{33}
\end{pmatrix}
\begin{pmatrix}
c_{11}&c_{12}&c_{13}\\
0&1&c_{23}\\
0&0&c_{33}
\end{pmatrix}
\begin{pmatrix}
y'\\v\\0
\end{pmatrix},
\label{eq:chi1}
\end{equation}
or, equivalently,
$$
\begin{pmatrix}
	d_{11}x'+d_{12}u\\u\\0
\end{pmatrix}
=
\begin{pmatrix}
	a_{11}&a_{12}&a_{13}\\
	a_{21}&a_{22}&a_{23}\\
	a_{31}&a_{32}&a_{33}
\end{pmatrix}
\begin{pmatrix}
	c_{11}y'+c_{12}v\\v\\0 
\end{pmatrix}.
$$
Thus for a given $u$, $v$ solutions $x$, $y$ and $x'$, $y'$ of systems
(\ref{eq:def-chi}) and (\ref{eq:chi1})
are connected by 
$$x=
d_{11}x'+d_{12}u
\qquad
y=c_{11}y'+c_{12}v.
$$
Since matrices $d_{11}$, $c_{11}$ are invertible, we get that $u$, $v$ in both cases
are same.


\sm

{\bf \punct The discrete invariant.}

\begin{proposition}
	\label{pr:aaa}
	 Numbers
$$
\rk \begin{pmatrix}
a_{31}
\end{pmatrix}, \qquad
\rk \begin{pmatrix}
a_{31}&a_{32}
\end{pmatrix}, \qquad  
\rk \begin{pmatrix}
a_{21}\\
a_{31}
\end{pmatrix},\qquad
\rk \begin{pmatrix}
a_{21}&a_{22}\\
a_{31}&a_{32}
\end{pmatrix}
$$
are invariants of double cosets.
\end{proposition}

{\sc Proof}
Indeed, let
$$ A'=
\begin{pmatrix}
d_{11}&d_{12}&d_{13}\\
0&1&d_{23}\\
0&0&d_{33}
\end{pmatrix}
\begin{pmatrix}
a_{11}&a_{12}&a_{13}\\
a_{21}&a_{22}&a_{23}\\
a_{31}&a_{32}&a_{33}
\end{pmatrix}
\begin{pmatrix}
c_{11}&c_{12}&c_{13}\\
0&1&c_{23}\\
0&0&c_{33}
\end{pmatrix}.
$$
Then
$$
a_{31}'= d_{33} a_{31} c_{11}, \qquad \begin{pmatrix}
a_{21}'&a_{22}'\\
a_{31}'&a_{32}'
\end{pmatrix}
=
\begin{pmatrix}
1&d_{23}\\0&d_{33}
\end{pmatrix}
\begin{pmatrix}
a_{21}&a_{22}\\
a_{31}&a_{32}
\end{pmatrix}
\begin{pmatrix}
c_{11}&c_{21}\\0&1
\end{pmatrix},
$$ 
etc. The statement becomes obvious.
\hfill $\square$

\sm

{\bf\punct Completeness of the system of invariants.}
Here we prove Theorem \ref{th:completeness}, i.e., show that {\it the characteristic
linear relation $\chi(\fra)$ and the invariant $\eta(\fra)\in \Z_+$
completely determine a double coset $\fra$.} 

Consider 'parabolic'  groups $P_\alphA\supset Q_\alphA$ consisting of matrices
$$
\begin{pmatrix}
c_{11}&c_{12}&c_{13}\\
0&c_{22}&c_{23}\\
0&0&c_{33}
\end{pmatrix}, \qquad c_{22}\in \GL(|\alphA|,\k).
$$
Clearly, $Q_\alphA$ is a normal subgroup in $P_\alphA$, the quotient is $\GL(|\alphA|,\k)$.
This implies the following observation:

\begin{lemma}
	Let $R$ ranges in $\GL(|\alphA|,\k)$, $S$ in $\GL(|\betA|,\k)$. Then the map
	$$A\mapsto \begin{pmatrix}
	1&0&0\\ 0&R&1\\0&0&1
	\end{pmatrix}^{-1} A \begin{pmatrix}
	1&0&0\\ 0&S&1\\0&0&1
	\end{pmatrix}$$
	induces an action of the group $\GL(|\alphA|,\k)\times \GL(|\betA|,\k)$
	on the double coset space $Q_\alphA\setminus \GL(2\infty,\k) /Q_\betA$.
\end{lemma} 

On the other hand the same group acts on the set of linear relations
$L:\k^{|\betA|}\rra \k^{|\alphA|}$ by
$$
L\mapsto R^{-1} L S.
$$

The following statement also is obvious.

\begin{lemma}
	The map $\fra \mapsto\chi(\fra)$ is $\GL(|\alphA|,\k)\times \GL(|\betA|,\k)$-equivariant.
\end{lemma}

Let us describe double cosets
$P_\alphA\setminus\GL(2\infty,\k)/P_\betA$.

\begin{lemma}
	\label{l:gauss}
	Any double coset in $P_\alphA\setminus \GL(2\infty,\k)/P_\betA$
	has a unique representative as a {\rm 0-1}-matrix of the form
	\begin{equation}
	{J_\kappa=
	\vphantom{\left(\begin{array}{ccc|ccc|ccc}
	1\\
	0\\
	0\\
	\hline
	0\\
	0\\
	0\\
	\hline
	0\\
	0\\
	0
	\end{array}\right)}}^\nw
	\left(\begin{array}{ccc|ccc|ccc}
	\1&0&0 &0&0&0 &0&0&0\\
	0&0&0 &\1&0&0 &0&0&0\\
	0&0&0 &0&0&0 &\1&0&0\\
	\hline
	0&\1&0 &0&0&0 &0&0&0\\
	0&0&0 &0&\1&0 &0&0&0\\
	0&0&0 &0&0&0 &0&\1&0\\
	\hline
	0&0&\1 &0&0&0 &0&0&0\\
	0&0&0 &0&0&\1 &0&0&0\\
	0&0&0 &0&0&0 &0&0&\1
	\end{array}\right)_\se
	,
	\label{eq:J-kappa}
	\end{equation}
	where sizes $\kappa_{ij}\ge 0$ of units in $ij$-blocks satisfy
	conditions
\refstepcounter{equation}
\label{eq:step}	
	\begin{align*}
\kappa_{11}+\kappa_{12}+\kappa_{13}&=M_-;
\label{eq:row1}
 \tag{\ref{eq:step}.row1}
\\
 \kappa_{21}+\kappa_{22}+\kappa_{23}&=|\alphA|=\alpha_+-\alpha_-;
 \label{eq:row2}
  \tag{\ref{eq:step}.row2}
 \\
\kappa_{31}+\kappa_{32}+\kappa_{33}&=M_+ ;
\label{eq:row3}
 \tag{\ref{eq:step}.row3}
\\
\kappa_{11}+\kappa_{21}+\kappa_{31}&=N_-;
\label{eq:col1}
 \tag{\ref{eq:step}.col1}
\\
 \kappa_{12}+\kappa_{22}+\kappa_{33}&=|\betA|=\beta_+-\beta_-;
 \label{eq:col2}
 \tag{\ref{eq:step}.col2}
\\
\kappa_{13}+\kappa_{23}+\kappa_{33}&=N_+
.
\label{eq:col3}
\tag{\ref{eq:step}.col3}
\end{align*}
\end{lemma}

Recall that $\alpha_\pm$, $\beta_\pm$ are fixed and $M_\pm$, $N_\pm$
satisfy conditions
\refstepcounter{equation}
\label{eq:step1}	
\begin{align}
M_--\alpha_-=N_--\beta_-;
\label{eq:MMNN1}\tag{\ref{eq:step1}$_-$}
\\
M_++\alpha_+=N_++\beta_+.
\label{eq:MMNN2}
\tag{\ref{eq:step1}$_+$}
\end{align}

The lemma follows from the Gauss reduction of systems of linear equation and we omit its proof.
\hfill $\square$

   
	
	\sm 
	

{\sc Remark.}
Replacing $\kappa_{11}\mapsto \kappa_{11}+1$ {\rm(}resp., $\kappa_{33}\mapsto \kappa_{33}+1${\rm)}
does not change the matrix $J_\kappa$ (due to the presence arrows $\nw$ , $\se$ in (\ref{eq:J-kappa})). In particular, 
we can set $\kappa_{11}=0$, $\kappa_{33}=0$.	
\hfill $\boxtimes$

\begin{lemma}
	The linear relation 
	$
	\chi(J_\kappa)
	$
	from
		$$
	\k_{\betA}=\k^{\kappa_{21}}\oplus\k^{\kappa_{22}}\oplus \k^{\kappa_{23}}
	$$
to	
	$$
	\k_{\alphA}=\k^{\kappa_{12}}\oplus\k^{\kappa_{22}}\oplus \k^{\kappa_{32}}
	$$	
	consists of all vectors of the form
	$$
	(v,u,0)\oplus (w,u,0).
	$$
	In particular,
	\begin{equation}
	\rk \chi(J_\kappa)=\kappa_{22},\quad \dim \indef \chi(J_\kappa)=\kappa_{21},
	\quad \dim\ker \chi(J_\kappa)=\kappa_{12}.
	\label{eq:rk}
	\end{equation}
\end{lemma}

This follows from  a straightforward calculation.
\hfill $\square$

\sm 

Notice that each orbit of the group $\GL(|\alphA|,\k)\times\GL(|\betA|,\k)$
on the set of linear relations $\k^\betA\rra \k^\alphA$ has a unique
representative of the form $\chi(J_\kappa)$.

Theorem \ref{th:completeness} is a corollary of the following lemma.

\begin{lemma}
	The map $\fra\mapsto\chi(\fra)$ is a bijection on each $\GL(|\alphA|,\k)\times\GL(|\betA|,\k)$-orbit.
\end{lemma}

{\sc Proof.} It is sufficient to show that  the stabilizer $\cM(J_\kappa)$ of 
a double coset $Q_\alphA\cdot J_\kappa\cdot Q_\betA$
coincides with the stabilizer $\cN(J_\kappa)$ of the linear relation  $\chi(J_\kappa)$. 
The inclusion $\cM(J_\kappa)\subset \cN(J_\kappa)$ follows from the equivariance.
Let us prove the inclusion inclusion.

The stabilizer $\cN(\kappa)$ consists
of pairs $(R,S)\in \GL(|\alphA|,\k)\times \GL(|\betA|,\k)$ having the form
$$
\begin{pmatrix} r_{11}&r_{12}&r_{13}\\
0&r_{22}&r_{23}\\
0&0&r_{33}
\end{pmatrix}, \qquad
\begin{pmatrix} s_{11}&s_{12}&s_{13}\\
0&s_{22}&s_{23}\\
0&0&s_{33}
\end{pmatrix},\qquad \text{where $r_{22}=s_{22}$.}
$$
Indeed,
 a matrix $S$ must preserve the flag $\ker \chi(J_\kappa)\subset \dom \chi(J_\kappa)$,
a matrix $R$ must preserve the flag $\indef \chi(J_\kappa)\subset \im \chi(J_\kappa)$.
This implies triangular forms of 
our matrices. The linear map
$$
\dom \chi(J_\kappa)/\ker \chi(J_\kappa)\to \im \chi(J_\kappa)/\indef \chi(J_\kappa)
$$
in our case is identical and this implies $r_{22}=s_{22}$.

Let us show that such pairs stabilize the double coset $Q_\alphA\cdot J_\kappa\cdot Q_\betA$.
Without loss of generality we can assume $\kappa_{11}=\kappa_{33}=0$.
Denote $T=R^{-1}$, so $t_{22}=s_{22}^{-1}$. Then
$$
\begin{pmatrix} 1&0&0\\
0&T&0\\
0&1&0
\end{pmatrix} \cdot J_\kappa \cdot \begin{pmatrix} 1&0&0\\
0&S&0\\
0&1&0
\end{pmatrix}
=\left(\begin{array}{cc|ccc|cc}
0&0&s_{11}&s_{12}&s_{13}&0&0\\
0&0&0&0&0&1&0\\
\hline 
t_{11}&0&0&t_{12}s_{22}&t_{12}s_{23}&0&t_{13}\\
0&0&0&1&t_{22}s_{23}&0&t_{23}\\
0&0&0&0&0&0&t_{33}\\
\hline
0&1&0&0&0&0&0\\
0&0&0&0&s_{33}&0&0
\end{array}\right).
$$
It is more or less clear that multiplying such matrices by elements
of $Q_\alphA$ from the left and elements of $Q_\betA$ from the right
we can reduce this matrix to the form $J_\kappa$.
Formally, the last product is equal to
{\small
\begin{multline*}
\left(
\begin{array}{cc|ccc|cc}
s_{11}&0&0&s_{12}&-s_{12}t_{23}t_{33}^{-1}&0&s_{13}\\
0&1&0&0&0&0&0\\
\hline 
0&0&1&0&0&0&0\\
0&0&0&1&0&0&t_{22}s_{23}\\
0&0&0&0&1&0&0\\
\hline
0&0&0&0&0&1&0\\
0&0&0&0&0&0&t_{33}
\end{array}
\right)\times \\ \times
J_\kappa
\left(
\begin{array}{cc|ccc|cc}
t_{11}&0&0&t_{12}s_{22}&t_{12}s_{23}&0&t_{13}\\
0&1&0&0&0&0&0\\
\hline
0&0&1&0&0&0&0\\
0&0&0&1&0&0&t_{23}\\
0&0&0&0&1&0&0\\
\hline 
0&0&0&0&0&1&0\\
0&0&0&0&0&0&t_{33}
\end{array}
\right)
\end{multline*}
}
and this completes the proof.
\hfill $\square$

\sm 

{\bf\punct The expression for $\eta(\fra^*)$.}
Here we prove the statement c) of Theorem \ref{th:isomorphism}.
%
%
%
%
We can assume that $A$ has  form $J_\kappa$, see (\ref{eq:J-kappa}).
Extracting $\beta_-$ from both sides of (\ref{eq:row1}) and $\alpha_-$ from both sides of
(\ref{eq:col1}) and keeping in the mind the equality $N_--\beta_-=M_--\alpha_-$,
see (\ref{eq:MMNN1}), we come to
$$
\kappa_{11}+\kappa_{21}+\kappa_{31}-\beta_-= \kappa_{11}+\kappa_{12}+\kappa_{13}-\alpha_-,
$$
or
\begin{align}
\eta(\fra^*)&=\kappa_{13}=\kappa_{21}+\kappa_{31}-\beta_- -\kappa_{12}+\alpha_-=
\label{eq:eta*1}
\\&=
\dim \indef \chi(\fra) +\eta(\fra)-\beta_- -\dim\ker \chi(\fra)+\alpha_+
\label{eq:eta*2}
\end{align}
The last line  is formula (\ref{eq:fra*}).

\sm

{\bf \punct Inequalities for $\eta(\fra)$.} Here we prove Proposition \ref{pr:equality}
about possible domain for $\eta(\fra)$ if $\chi(\fra)$ is fixed.
The expression in the line (\ref{eq:eta*2}) must be $\ge 0$, and this
implies the desired inequality (\ref{eq:chi-eta}). We must show that
 this is sufficient.
 
 \begin{lemma}
 	\label{eq:13}
 	Denote by $\Xi$ the set of all $13$-ples of
 	$$  
 	\kappa_{ij}\in \Z_+, \quad\text{where $1\le i,j\le 3$ and}\qquad \text{$N_\pm$, $ M_\pm\in \Z_+$}.
 	$$
 satisfying $10$ equations {\rm (\ref{eq:row1})--(\ref{eq:col3}), (\ref{eq:MMNN1})--(\ref{eq:MMNN2})}. Then 
 possible sub-tuples $(\kappa_{21}, \kappa_{22}, \kappa_{31}, \kappa_{12})$
 are precisely  integer points  of the cone $\Delta$ defined by equalities
 \refstepcounter{equation}
 \label{eq:step3}
  \begin{align}
  \kappa_{21}\ge 0,\quad \kappa_{22}\ge 0,\quad \kappa_{31}\ge 0,\quad \kappa_{12}\ge 0;
   \label{eq:Delta1}
   \tag{\ref{eq:step3}.ineq1}
  \\
  \kappa_{12}+\kappa_{22}\le \beta_+-\beta_-,\qquad  \kappa_{12}+\kappa_{22}\le  \alpha_+-\alpha_-;
  \label{eq:Delta2}
  \tag{\ref{eq:step3}.ineq2}
  \\
  \kappa_{21}+\kappa_{31}-\beta_- -\kappa_{12}+\alpha_-\ge 0.
  \label{eq:Delta3}
  \tag{\ref{eq:step3}.ineq3}
  \end{align} 
 \end{lemma}

{\sc Proof.} All steps of the proof are obvious but the result is not clear
until steps are  preformed. 

We notice that our 8 equations are dependent:
the sum of 3 equations (\ref{eq:row1})--(\ref{eq:row3}) minus the sum of 3 equations
(\ref{eq:col1})--(\ref{eq:col2}) coincides with the sum of (\ref{eq:MMNN1}) and (\ref{eq:MMNN2}).

Next, fix a point  $(\kappa_{21}, \kappa_{22}, \kappa_{31}, \kappa_{12})\in\Delta$
 and construct a point of $\Xi$ over it. Let assign the remaining coordinates step by step.

\sm
 
1) $\kappa_{23}$, $\kappa_{32}$. We find them from the equations (\ref{eq:col2}) 
and (\ref{eq:row2}). By (\ref{eq:Delta2}), $\kappa_{23}$, $\kappa_{32}\in \Z_+$.

\sm 

2) Set $\kappa_{11}=\kappa_{33}=0$.

\sm 

3) $M_+$, $N_-$. We find them from  equations (\ref{eq:row3}) and (\ref{eq:col1}).
Obviously, they are in $\Z_+$.

\sm 

4) $M_-$. We evaluate it from the equation (\ref{eq:MMNN1}). Positivity of $M_-$ in this moment
is not obvious.

\sm 

5) We find $\kappa_{13}$ from the equation (\ref{eq:row1}) and  get $\kappa_{13}=\kappa_{21}+\kappa_{31}-\beta_- -\kappa_{12}+\alpha_-$ (in fact this calculution
is present in the previous subsection). The condition (\ref{eq:Delta3}) claims that it is positive.
Therefore $M_-$ is positive by (\ref{eq:row1}). 

\sm 

6) $N_+$. We find it from (\ref{eq:col3}). Obviously, $N_+\in \Z_+$.

\sm

Thus we get a vector in $\Z_+^{13}$. We  used 7 equations, (\ref{eq:row1})--(\ref{eq:col3})
and (\ref{eq:MMNN1}). So they are satisfied. The 8-th equation is satisfied automatically.
\hfill $\square$

\sm

{\bf\punct Characteristic linear relations of  products of double cosets.}
Let 
\begin{equation}
\begin{pmatrix}
x_2\\u\\0
\end{pmatrix}
=
\begin{pmatrix}
a_{11}&a_{12}&a_{13}\\
a_{21}&a_{22}&a_{23}\\
a_{31}&a_{32}&a_{33}
\end{pmatrix}
\begin{pmatrix}
y_2\\v\\0
\end{pmatrix},
\quad\begin{pmatrix}
x_1\\v\\0
\end{pmatrix}
=
\begin{pmatrix}
b_{11}&b_{12}&b_{13}\\
b_{21}&b_{22}&b_{23}\\
b_{31}&b_{32}&d_{33}
\end{pmatrix}
\begin{pmatrix}
y_1\\w\\0
\end{pmatrix}.
\label{eq:AB}
\end{equation}
Then
\begin{equation}
A^\circ B^\lozenge \begin{pmatrix}
y_1\\y_2\\w\\0\\0
\end{pmatrix}
=
A^\circ 
\left( \begin{array}{ccccc}
b_{11}&0&b_{12}&0&b_{13}\\
0&1&0&0&0\\
b_{21}&0&b_{22}&0&b_{23}\\
0&0&0&1&0\\
b_{31}&0&b_{32}&0&b_{33}
\end{array}
\right)
\begin{pmatrix}
y_1\\y_2\\w\\0\\0
\end{pmatrix}
=A^\circ  \begin{pmatrix}
x_1\\y_2\\v\\0\\0
\end{pmatrix}=
\label{eq:high1}
\end{equation}
\begin{equation}
\qquad
=
\left(\begin{array}{ccccc}
1&0&0&0&0\\
0&a_{11}&a_{12}&a_{13}&0\\
0&a_{21}&a_{22}&a_{23}&0\\
0&a_{31}&a_{32}&a_{33}&0\\
0&0&0&0&1
\end{array}\right)
\begin{pmatrix}
x_1\\y_2\\v\\0\\0
\end{pmatrix}
= \begin{pmatrix}
x_1\\x_2\\u\\0\\0
\end{pmatrix}.
\label{eq:high2}
\end{equation}

Therefore $(w,u)\in \chi(\fra \star \frb)$. 

\sm 

Conversely, let the right hand side of the equation (\ref{eq:high1})--(\ref{eq:high2})
equals to the left hand side. Then applying $B^\lozenge$ to a column
$\begin{pmatrix}
y_1&y_2&w&0&0
\end{pmatrix}^t$, we get a column of the form
$\begin{pmatrix}
z_1&y_2&q&0&s
\end{pmatrix}^t$. Applying $A^\circ$ to this column, we get an expression of the form
$\begin{pmatrix}
z_1&z_2&r&t&s
\end{pmatrix}^t$. But we must get a vector of a form $\begin{pmatrix}
x_1&x_2&u&0&0
\end{pmatrix}^t$, i.e., $t=0$, $s=0$.
Hence
\begin{equation*}
\begin{pmatrix}
z_2\\r\\0
\end{pmatrix}
=
\begin{pmatrix}
a_{11}&a_{12}&a_{13}\\
a_{21}&a_{22}&a_{23}\\
a_{31}&a_{32}&a_{33}
\end{pmatrix}
\begin{pmatrix}
y_2\\q\\0
\end{pmatrix},
\quad\begin{pmatrix}
x_1\\q\\0
\end{pmatrix}
=
\begin{pmatrix}
b_{11}&b_{12}&b_{13}\\
b_{21}&b_{22}&b_{23}\\
b_{31}&b_{32}&d_{33}
\end{pmatrix}
\begin{pmatrix}
y_1\\w\\0
\end{pmatrix}.
\end{equation*}
Therefore
$(w,q)\in \chi(\frb)$, $(q,r)\in \chi(\fra)$.
\hfill $\square$

\sm 

{\bf\punct The invariant $\eta$ of a product of double cosets.}
It remains to prove the formula
\begin{equation}
\eta(\fra\star\frb)=\eta(\fra)+\eta(\frb)+\dim \indef \frb/ (\indef \frb\cap\dom\fra ).
\label{eq:eta-ab}
\end{equation}
Let us pass to another invariant
$$
\xi(\fra):=\rk\begin{pmatrix}
a_{21}\\ a_{31}
\end{pmatrix}=\eta(\fra)+\dim \indef\chi(\fra)
$$
(the identity is clear from the canonical forms).
It is easy to see that for any linear relations
$P:X\rra Y$, $Q:Y\rra Z$ we have
$$
\dim \indef QP=\dim \indef Q+ \dim (\indef P\cap \dom Q)-\dim (\indef P\cap \ker Q).
$$
Hence
\begin{multline*}
\dim \indef\bigl(\chi(\fra)\chi(\frb)\bigr)
=\dim \indef \chi (\fra)+\dim\indef \chi(\frb) -\\
-\dim \bigl( \indef \chi(\frb)/(\indef \chi(\frb)\cap \dom \chi(\fra))\bigr)
-\dim \bigl(\indef \chi(\frb)\cap \ker\chi(\fra )\bigr).
\end{multline*}
Therefore (\ref{eq:eta-ab}) can be written as
\begin{equation}
\xi(\fra\star\frb)=\xi(\fra)+\xi(\frb)-\dim\bigl(\ker \chi(\fra)\cap \indef \chi(\frb)\bigr).
\label{eq:xi-ab}
\end{equation}
We wish to prove the last identity.

\sm 

Consider the space $\cM(\fra)$ of all $y$, for which there exists $x$ such that
\begin{equation}
\begin{pmatrix}
x\\0\\0
\end{pmatrix}
=
\begin{pmatrix}
a_{11}&a_{12}&a_{13}\\
a_{21}&a_{22}&a_{23}\\
a_{31}&a_{32}&a_{33}
\end{pmatrix}
\begin{pmatrix}
y\\0\\0
\end{pmatrix}.
\label{eq:def-xi}
\end{equation}
We can say the same in a shorter way,
$\cM$ is defined by the equation
$$
0=
\begin{pmatrix}
a_{21}\\a_{31}
\end{pmatrix}
\begin{pmatrix}
y
\end{pmatrix}.
$$
Clearly, $\xi(\fra)$ is the codimension of $\cM$ in the space of all $y$.

We must evaluate the codimension of the subspace $\cM(\fra \star\frb)$
 of all $(y_1, y_2)$ such that 
 there exists $(x_1, x_2)$ satisfying
 \begin{equation*}
 \begin{pmatrix}
 x_1\\x_2\\0\\0\\0
 \end{pmatrix}
 =
 \begin{pmatrix}
 A^\circ B^\lozenge 
 \end{pmatrix}
 \begin{pmatrix}
 y_1\\y_2\\0\\0\\0
 \end{pmatrix}.
 \end{equation*}
 Applying the matrix $B^\lozenge$ to a vector $\begin{pmatrix}
 y_1&y_2&0&0&0
 \end{pmatrix}^t$ we get a vector of the form $\begin{pmatrix}
 z_1&y_2&p&0&h_2
 \end{pmatrix}^t$. Applying $A^\circ$ we come to a vector of the form
 $\begin{pmatrix}
 z_1&z_2&q&h_1&h_2
 \end{pmatrix}^t$. We want $h_1=0$, $h_2=0$, and $q=0$. Therefore we have
 \begin{equation}
 \begin{pmatrix}
 x_2\\0\\0
 \end{pmatrix}=\begin{pmatrix}
 a_{11}&a_{12}&a_{13}\\
 a_{21}&a_{22}&a_{23}\\
 a_{31}&a_{32}&a_{33}
 \end{pmatrix}\begin{pmatrix}
 y_2\\p\\0
 \end{pmatrix}, \quad
 \begin{pmatrix}
 x_1\\p\\0
 \end{pmatrix}=
 \begin{pmatrix}
 b_{11}&b_{12}&b_{13}\\
 b_{21}&b_{22}&b_{23}\\
 b_{31}&b_{32}&b_{33}
 \end{pmatrix}
 \begin{pmatrix}
 y_1\\0\\0
 \end{pmatrix}.
 \label{eq:for-xi}
 \end{equation}
We see that $p\in \indef \chi(\frb)\cap \ker \chi(\fra)$.

Clearly $\cM(\fra \star \frb)\supset \cM(\fra)\oplus \cM(\frb)$.
Moreover, we have a surjective map
$$
\pi: \cM(\fra \star \frb) \to
 \indef \chi(\frb)\cap \ker \chi(\fra),
$$
and $\ker \pi\supset \cM(\fra)\oplus \cM(\frb)$.
Conversely, let $\begin{pmatrix} y_1\\y_2
\end{pmatrix}
\in \ker \pi$.
Then it satisfies two equations (\ref{eq:for-xi})
with $p=0$. Therefore $y_1\in \cM(\frb)$ and $y_2\in \cM(\fra)$.
and
$$
\cM(\fra \star \frb)\bigl/ \bigl(\cM(\fra)\oplus \cM(\frb)\bigr)\simeq
 \indef \chi(\frb)\cap \ker \chi(\fra).
$$
This completes the proof of Theorem \ref{th:isomorphism}.b.

\section{The group $\GLL$ and multiplicativity\label{s:4}}

\COUNTERS

Here we prove statements of Subsect. \ref{ss:i-multiplicativity}.
The key place is 'Mautner phenomenon', see Subsect.  \ref{ss:mautner}
and Lemma \ref{l:two-subspaces}. After this the proof of Theorem \ref{th:multiplicativity}
(multiplicativity)
becomes automatic.

Subsection \ref{ss:digression} contains  two observations outside the main topic of the paper,
the first is devoted to the group $\GLL(V\sqcup V^\diamond)$ (see Subsect. \ref{ss:review}),
 the second 
is related to groups of infinite matrices over $p$-adic integers.

\sm 

{\bf \punct Proof of Lemma \ref{l:repr-finite}.}

\sm 

a) 
We must  show that {\it for any $s\in \GLL$ the double coset
$\ov Q_\alphA\cdot s\cdot \ov Q_\betA$ contains a finitary matrix}. Without loss
of generality we can assume that $\alphA=\betA$. Otherwise we take $\gammA$ such that $\gammA\succ\alphA$,
$\gammA\succ\betA$
and  examine the double coset $\ov Q_\gammA\cdot s\cdot \ov Q_\gammA$.
Thus let represent $s$ as a block matrix of size $(\infty+|\alphA|+\infty)\times(\infty+|\alphA|+\infty)$,
$$
s=\begin{pmatrix}
a&b&c\\d&e&f\\g&h&k
\end{pmatrix}.
$$
Transformations
$$
s\mapsto 
\begin{pmatrix} u_1&0&0\\0&1&0\\0&0&v_1\end{pmatrix}^{-1} s
\begin{pmatrix} u_2&0&0\\0&1&0\\0&0&v_2\end{pmatrix}
$$
send
$$
s\mapsto u_1^{-1} a u_2, \qquad k\mapsto k_1^{-1} k v_2.
$$
The matrices $a$, $k$ are Fredholm matrices in the sense
of \cite{Ner-finite}, Subsects. 2.4--2.7, their Fredholm indices are 0.
Therefore we can reduce $a$ and $k$ to the forms
$$
a=\begin{pmatrix}1_\infty&0\\0&\0\end{pmatrix},
\qquad d=\begin{pmatrix}
\0&0\\0&1_\infty
\end{pmatrix},
$$
where two $\0$'s are square matrices  (see \cite{Ner-finite}, Lemma 2.7).
Hence our double coset contains a matrix of the form
$$
s'=\begin{pmatrix}
 1&0&b_1&c_{11}&c_{12}\\
  0&0&b_2&c_{21}&c_{22}\\
  d_1&d_2&e&f_1&f_2\\
  g_{11}&g_{12}&h_{1}&0&0\\
  g_{21}&g_{22}&h_{1}&0&1
\end{pmatrix}.
$$
Multiplying such matrices from the left and right by matrices of the form
$$
\begin{pmatrix}
1&0&*&*&*\\
0&1&0&0&*\\
0&0&1&0&*\\
0&0&0&1&0\\
0&0&0&0&1
\end{pmatrix}\in \ov Q_\alphA, 
$$
we can make zero from $b_1$, $c_{11}$, $c_{12}$, $c_{22}$, $f_2$.
The blocks $e$, $b_2$, $c_{21}$, $f_1$ have finite sizes, 
the definition of $\GLL$  implies that the blocks
$\begin{pmatrix} d_1&d_2 \end{pmatrix}$, $\begin{pmatrix}
g_{11}&g_{12}\\
g_{21}&g_{22}\end{pmatrix}$, $\begin{pmatrix}
e_1\\e_2
\end{pmatrix}$
contain only  finite numbers of nonzero matrix elements.
Thus we get a finitary matrix.

\sm 

b) By Theorem \ref{th:completeness}, the invariants $\chi(\fra)$ and $\eta(\fra)$ separate
double cosets $Q_\alphA\setminus \GL(2\infty,\F)/Q_\betA$.
However, $\chi(\fra)$ and $\eta(\fra)$ also are invariants of double cosets 
$\ov Q_\alphA\setminus \GLL/\,\ov Q_\betA$ (our proof in Subsect. \ref{ss:invariants-proofs}
is valid in this case).

\sm 

{\bf \punct Proof of Lemma \ref{l:admissibility}.%
\label{ss:adimssibility-proof}}
Recall that $H_\alphA$ denotes the space of $Q_\alphA$-fixed vectors
in a representation of $\GL(\infty,\F)$. We must show that
{\it for a unitary representation of $\GL(\infty,\F)$ the continuity
in the topology of $\GLL$ is equivalent to density of the space
$\cup_{\alphA}H_\alphA$.}  

\sm 

{\sc The statement $\Rightarrow$.}
The subgroups $\ov Q_\alphA\subset \GLL$ are open and form a fundamental system
of neighborhoods of unit in $\GLL$. This is sufficient for application
of Proposition VIII.1.2 from \cite{Ner-book}, which immediately gives the desired 
statement.

\sm

{\sc  The statement $\Leftarrow$.}
Conversely, let $\cup_{\alphA}H_\alphA$ be dense. We must verify that
matrix elements $\la\rho(g) h_1, h_2\ra$ are continuous in the topology of $\GLL$.
It is sufficient to do this for $h_1$, $h_2$ ranging in a dense subspace
in $H$, in particular in the subspace $\cup H_\alphA$. However,
if $h_1\in H_\betA$, $h_2\in H_\gammA$, then our matrix element
is a function on a countable space 
$$
Q_\gammA\setminus \GL(2\infty,\F)/Q_\betA\simeq \ov Q_\gammA\setminus \GLL/\,\ov Q_\betA.
$$ 
Since subgroups $\ov Q_\deltA\subset \GLL$ are open, the double coset space in the right hand
side is discrete, and all functions on this space are continuous.

\sm 

{\bf\punct The Mautner phenomenon. Coincidence of spaces of fixed vectors.%
\label{ss:mautner}}
Recall the following phenomenon related to Lie groups, which was discovered by Gelfand and Fomin
\cite{GF}
and investigated in details by Mautner and Moore, see \cite{Moo}. Let $G$ be a Lie group, $H$ a non-compact subgroup.
Then very often a vector in a unitary representation fixed by $H$ is automatically fixed  by
a larger subgroup  $\wh H\subset G$.

Recall that $\GLL(V)$ denotes the group of all linear operators 
in the countable linear space $V$ over $\F$. By $\GLL(V^\diamond)$ we denote
the group of all continuous linear operators in the dual space $V^\circ$, see Subsect. \ref{ss:review}.
 Both groups are present in $\ov Q_\alphA$ as subgroups consisting of matrices of the form
 $$
 \begin{pmatrix}
 1&0&0\\0&1_{|\alphA|}&0\\0&0&d
 \end{pmatrix}\qquad \text{and} \qquad
  \begin{pmatrix}
 a&0&0\\0&1_{|\alphA|}&0\\0&0&1
 \end{pmatrix}
 $$
 respectively.

\begin{lemma}
	\label{l:SGL}
	For a unitary representation of the group $\GLL(V)$, any
	$S_\infty$-fixed vector is fixed by the whole group $\GLL(V)$. 
\end{lemma}

The statement can be  derived from Tsankov's
 classification \cite{Tsa} of unitary representations of $\GLL(V)$,
 however, we  present a simple direct proof. 

For each group $\ov Q_\alphA$ consider the subgroup $\ov S_{[\alphA]}$ consisting
of 0-1-matrices of the form
$$
\begin{pmatrix}
*&0&0\\
0&1_{|\alphA|}&0\\
0&0&*
\end{pmatrix}
.$$
By $S_{[\alphA]}$ denote its (dense) subgroup consisting of finitary matrices.

\begin{lemma}
\label{l:two-subspaces}
	Let $\rho$ be a unitary representation of the group
	$\ov Q_\alphA$ in a Hilbert space $H$. Let $h\in H$ be an $S_{[\alphA]}$-fixed vector.
	Then $h$ is $\ov Q_{\alphA}$-fixed.
\end{lemma}

Proofs of these lemmas are based on the following statement.

\begin{proposition}
	\label{pr:mautner}
	Let a countable discrete group $\Gamma$ act by automorphisms on a compact  Abelian group
	 $N$. Let $\wh N$ be the Pontryagin dual group, i.e. the group of characters
	 of $N$. Assume that all orbits of $\Gamma$ on the discrete group
	  $\wh N$ except the orbit of the trivial character are infinite. Then for any unitary representation
	 of the semidirect product $\Gamma\ltimes N$, any $\Gamma$-fixed vector is  fixed by the whole group $\Gamma\ltimes N$.
\end{proposition}

{\sc Proof of Proposition \ref{pr:mautner}.}  The group $\Gamma\ltimes N$ is locally compact. Therefore a unitary representation of this group
can be decomposed into a direct integral of irreducible representations (see, e. g., \cite{Mac}, Sect 2.6,
\cite{Kir}, Subsect. 8.4).
We claim that any irreducible representation $\rho$ of $\Gamma\ltimes N$ having an $\Gamma$-fixed vector
is trivial. 

According the Mackey theorem about a unitary dual of locally compact group with an Abelian normal subgroup (see, e.g., \cite{Kir}, Theorem 13.3.1),
any irreducible unitary representations of the group $\Gamma\ltimes N$ can be realized in the following way.
Consider an
orbit $\Omega$ of $\Gamma$ on $\wh N$, fix $\chi_0\in \Omega$.  Denote by $\Delta$ the stabilizer
of $\chi_0$ in $\Gamma$, fix an irreducible unitary representation
$\tau$ of $\Delta$
in a Hilbert space $K$.
Consider the space $\ell^2(\Omega,K)$ of 
$\ell^2$-functions on the discrete set $\Omega=\Delta\setminus\Gamma$ taking values in $K$.
The Abelian
subgroup $N$ acts in this space by multiplications 
\begin{equation} n: \, F(\chi)\mapsto \chi(n) F(\gamma) ,
\label{eq:chiF}
\end{equation}
where $\chi(n)$ denotes the value of a character $\chi\in \wh N$ on
an element $n\in N$.    
The group $\Gamma$ acts by transformations
of the form
$$
\gamma: F(\chi)\mapsto T(\gamma,\chi)
F(\chi\gamma),
$$ 
where $T$ is a function from $\Gamma\times \Omega$
to the unitary group of the space $K$ satisfying the cocycle identity
\begin{equation*}
T(\gamma_1\gamma_2,\chi)=T\bigl(\gamma_1, \chi \gamma_2)
\, T(\gamma_2, \chi)
\end{equation*}
and the condition
$$
T(\gamma,\chi_0)=\tau(\gamma)\qquad \text{for $\gamma\in\Delta$.}
$$

The norm of a function $F\in \ell^2(\Omega,K)$ is given by
$$
\|F\|_{\ell^2(\Omega,K)}^2=\sum_{\chi\in\Omega}\|F(\chi)\|^2_K.
$$
If $F$ is $\Gamma$-fixed, then all summands in the right-hand side coincide.
Therefore $F=0$ or $\Omega$ consists of one point (the trivial character),
$\ell^2(\Omega,K)$ is $K$.
By (\ref{eq:chiF}), the representation is trivial on the normal divisor $N$.
\hfill $\square$

\sm 

{\sc Proof of Lemma \ref{l:SGL}.}
Consider the subgroup in $\GLL(V)$ generated by $S_\infty$ and the group $N$ of
all diagonal matrices, so $N$ is a countable direct product of multiplicative groups
$\F^\times \simeq \Z_{q-1}$. Applying Proposition \ref{pr:mautner} to the  group $S_\infty\ltimes (\F^\times)^\infty$, we get that a vector fixed by 
$S_\infty$ is also fixed by all diagonal matrices.

Next, we consider the subgroup $Z\subset \GLL(V)$, consisting of all
  block  $(1+\infty)$-matrices having the form
 $
 \begin{pmatrix}
 1&x\\0&\sigma
 \end{pmatrix}
$,
   where $\sigma$ ranges in the group of finitary 0-1-matrices,
 and $x=\begin{pmatrix}
 x_1&x_2&\dots
 \end{pmatrix}$
 is arbitrary. So $x$ is contained in the direct product of a countable number of copies of $\F$. Applying Proposition \ref{pr:mautner} to this group, we get
 that a vector fixed by $S_\infty$ is also fixed by all matrices of the form
 $ \begin{pmatrix}
 1&x\\0&1
 \end{pmatrix}$. In particular, we can choose $x=(s,0, 0\dots)$. Conjugating
 this matrix by elements of $S_\infty$ we can get arbitrary matrices
 of the form $1+sE_{kl}$, where $k\ne l$ and $E_{kl}$ is the matrix having 1 on $kl$-th place. 
 
 Therefore a vector fixed by $S_\infty$ is fixed by all Chevalley generators of $\GL$. Hence it is fixed by the whole group $\GL(\infty,\F)$. By continuity,
 it is fixed by $\GLL(V)$.
 \hfill $\square$
 
 \sm
 
 {\sc Proof of Lemma \ref{l:two-subspaces}.} We apply Proposition \ref{pr:mautner}
 to two subgroups $H_1$, $H_2$ consisting of matrices
 $$
 \begin{pmatrix}
 \sigma_1 &*&*\\
 0&1&0\\
 0&0&1
 \end{pmatrix} \qquad\text{and}\qquad
  \begin{pmatrix}
 1 &0&*\\
 0&1&*\\
 0&0&\sigma_2
 \end{pmatrix}
 $$
 respectively, where $\sigma_1$, $\sigma_2$ a finitary 0-1-matrices. This implies that a $S_\alphA$-fixed vector $\xi$ is fixed by the subgroups consisting of all matrices
 of the form
  $$
   \begin{pmatrix}
 1 &*&*\\
 0&1&0\\
 0&0&0
 \end{pmatrix} \qquad\text{and}\qquad
    \begin{pmatrix}
 1 &0&*\\
 0&1&*\\
 0&0&1
 \end{pmatrix}
 $$
  respectively. Therefore $\xi$ is fixed by the product of these subgroups, i.e., by the group of all strictly upper triangular block matrices.
  
  It remains to apply Lemma \ref{l:SGL} to two subgroups $S_\infty$ in
  $Q_\alphA$.
  \hfill $\square$

  \sm
  
 {\bf \punct Some digressions. Admissibility in the Olshanski sense for the group 
 	$\GLL(V\sqcup V^\diamond)$ and its $p$-adic analogs.%
 \label{ss:digression}}

\sm 

{\sc The group $\GLL(V\sqcup V^\diamond)$.}
 Here we discuss some corollaries of Proposition \ref{pr:mautner} outside the main topic of this paper.
 Recall that $\GLL(V\sqcup V^\diamond)$ denotes the group of all infinite matrices over $\F$
 having a finite number of elements in each row and in each column, see Subsect. \ref{ss:review}.
 
 \begin{lemma}
 	\label{l:GLVV}
 	For  any unitary representation of the group $\GLL(V\sqcup V^\diamond)$
 	 any $S_\infty$-fixed vector is fixed by the whole group
 	$\GLL(V\sqcup V^\diamond)$.
 \end{lemma}
 
 {\sc Proof.} We can not literally repeat the
  proof
  the similar statement for $\GLL(V)$, i.e., Lemma \ref{l:SGL},
  since 
 the subgroup of $(1+\infty)$-block matrices $\begin{pmatrix}
 1&x\\0&1
 \end{pmatrix}$ now is not compact.
 We modify this place of the proof in the following way.
 
 Let us split the space $V$ as a direct sum of two-dimensional 
 subspaces, $V=\oplus_j W_j$. We regard subspaces $W_j$ as canonically isomorphic.
 Consider the subgroup $\Sigma\subset \GLL(V)$ consisting of finitary permutations
 of subspaces $W_j$. Consider the subgroup $\Delta$ consisting of block diagonal matrices, whose diagonal entries have the form $\begin{pmatrix}
 1&x_j\\0&1
 \end{pmatrix}$. We apply Proposition \ref{pr:mautner} to the semidirect product
 $\Sigma\ltimes \Delta$ and observe  that $\Delta$ also is contained in the stabilizer of $\xi$. Next, we set $x_1=s$, $x_2=x_3=\dots=0$ and get that 
 a Chevalley generator $1+sE_{12}$ also is contained in the stabilizer.
 The remaining part of the proof is the same. \hfill $\square$
 
 \sm 
 
 \begin{corollary}
 Unitary representations of $\GLL(V\sqcup V^\diamond)$ are admissible in the Olshanski sense.	
 \end{corollary}

{\sc Proof.}
 Denote by $\ov G^\alpha$ the subgroup in $\GLL(V\sqcup V^\diamond)$ consisting of block matrices of size
 $\alpha+\infty$ having the form $\begin{pmatrix}
 1&0\\0&*
 \end{pmatrix}$.
 Denote $\ov S_\infty^\alpha := \ov G^\alpha\cap \ov S_\infty$.
 For a unitary
 representation of $\GLL(V\sqcup V^\diamond)$ in a Hilbert space $H$ denote
 by $H[\alpha]$ the subspace of $\ov S_\infty^\alpha$-fixed vectors.
 Then $\cup_\alpha H[\alpha]$ is dense in $H$, see, e.g. \cite{Ner-book}, 
 Proposition VII.1.2. 
 It remains to  apply Lemma \ref{l:GLVV} to each subgroup
 $\ov G^\alpha$.
 \hfill $\square$
 
 \begin{lemma}
 Any double coset $\ov G_\alpha\setminus \GLL(V\sqcup V^\diamond)/G_\alpha$ contains
 a finitary matrix.	
 \end{lemma}	
 
 We omit a non-interesting proof (in particular the statement follows from
 slightly more difficult $p$-adic Lemma 4.1.a from \cite{Ner-p-adic}).
 
 \begin{corollary}
 	A unitary representation of $\GL(\infty,\F)$ admissible in the Olshanski sense
 	admits a continuous extension to $\GLL(V\sqcup V^\diamond)$.
 \end{corollary}
 
 It is sufficient to apply the argument from
Subsect. \ref{ss:adimssibility-proof}, the second statement.
\hfill $\square$

 \sm  
 
 {\sc Groups of infinite $p$-adic matrices with integer elements.}
 Now let $\frr$ be a compact commutative local ring. We keep in mind rings of integers 
 in locally compact non-Archimedean fields and their finite quotients as $\Z/p^n\Z$
 or truncated polynomial rings $\F[t]/t^n \F[t]$ (it seems that matrix groups over general local rings
 are not a topic of theory of unitary representations).
 Denote by $\frl=\frl(\frr)$ the module of all sequences in $\frr$ converging to 0.
 By $\GLL(\frl\sqcup\frl^\diamond)$
 we denote the group of matrices $g$ such that $g$, $g^t\in \GLL(\frl)$.
 
 \begin{lemma}
 	For any unitary representation of the group $\GLL(\frl\sqcup\frl^\diamond)$
 	any $S_\infty$-fixed vector is fixed by the whole group $\GLL(\frl\sqcup\frl^\diamond)$.
 \end{lemma} 

The proof of Lemma \ref{l:GLVV} remains to be valid in this case.

\sm

The lemma implies the Olshanski admissibility of unitary representations of
the group   $\GLL(\frl\sqcup\frl^\diamond)$. This statement is the
main result of the paper \cite{Ner-p-adic}.  
 
 \sm


\sm 

{\bf\punct Multiplicativity.}
Denote by $\ov S_\beta(\infty)\subset \ov S(\infty)$ the subgroup consisting of permutations
fixing $1$, \dots, $\beta\in\N$.
Denote
$$
I_N^{(\beta)}:=\begin{pmatrix}
          1_\beta&0&0\\
          0&0&1_N\\
          0&1_N&0
         \end{pmatrix}_\se
\in \ov S_\beta(\infty)
$$
We use the following statement (see \cite{Ner-book}, Theorem 1.4.c.

\sm

{\it Let $\rho$ be a unitary representation of the group $\ov S(\infty)$ in a Hilbert space $H$. 
Denote by $H_\beta\subset H$ the subspace of all  $\ov S_\beta(\infty)$-fixed vectors,
let $\Pi_\beta$ be the operator of orthogonal projection to $H_\beta$. Then
$\rho(I_N^{(\beta)})$ weakly converges to $\Pi_\beta$.}

\sm

The group $\ov S(\infty)$ has type $I$ (see \cite{Lie}), therefore $S(\infty)\times S(\infty)$
also has type $I$ (see, e.g., \cite{Mac}, Sect. 3.1). Therefore
(see \cite{Dix}, 13.1.8) irreducible unitary representations of $S(\infty)\times S(\infty)$
are tensor products of irreducible representations of factors. This implies the following
statement:

\begin{corollary}
  Let $\tau$ be a unitary representation of the group $\ov S(\infty)\times \ov S(\infty)$ in a Hilbert space $K$. Denote by $K_\beta\subset K$ the subspace of all  $\ov S_\beta(\infty)\times \ov S_\beta(\infty)$-fixed vectors,
let $\Pi_\beta$ be the operator of orthogonal projection to $K_\beta$. Then
$\tau(I_N^{(\beta)},I_N^{(\beta)})$ weakly converges to $\Pi_\beta$.
\end{corollary}

Let $J_\betA(\mu,\nu)\in \GLL$ be as above {\rm (\ref{eq:Jm})},
$$
J_\betA(\nu,\mu)=
\left(
\begin{array}{cc|c|cc}
0&1_\nu&0&0&0\\
1_\nu&0&0&0&0\\
\hline
0&0&1_{\betA}&0&0\\
\hline
0&0&0&0&1_\mu\\
0&0&0&1_\mu&0
\end{array}\right)
.$$

\begin{corollary}
	\label{cor:Jnn}
  Let $\rho$ be a unitary
 representation of the group $\GLL$ in a Hilbert space $H$. Then
 the sequence $\rho(J_\betA(n,n))$ weakly converges to $P_\betA$ as $n\to\infty$. 
\end{corollary}

{\sc Proof.} By the previous corollary
this sequence converges to the projector to the subspace of $\ov S_{[\betA]}$-fixed vectors.
By Lemma \ref{l:two-subspaces}, subspaces of $\ov S_{[\betA]}$-fixed vectors
and $\ov Q_\betA$-fixed vectors coincide.
\hfill $\square$

\sm

{\sc  Proof of Theorem \ref{th:multiplicativity}.}
Consider two double cosets $\fra\in \ov Q_\alphA\,\setminus \GLL/Q_\betA$,
$\frb\in \ov Q_\betA\,\setminus \GLL/Q_\gammA$. Choose finitary representatives
$A\in\fra$, $B\in\frb$ (see Lemma \ref{l:repr-finite}).
We must evaluate
\begin{multline*}
P_\alphA\, \rho(A)\, P_\betA \,\rho(B)\, P_\gammA=
\lim_{n\to\infty}
P_\alphA\, \rho(A)\, \lim_{n\to \infty} \rho\bigl(J_\beta(n,n)\bigr)\, \rho(B)\, P_\gammA=
\\
=
P_\alphA\,  \lim_{n\to \infty} \rho\bigl(A\,J_\beta(n,n)\,B\bigr)\, P_\gammA
\end{multline*}
This sequence is eventually constant. Its limit is 
$$
P_\alphA\,  \rho\bigl(A\,J_\beta(N,N)\,B\bigr)\, P_\gammA
$$
for sufficiently large $N$. By (\ref{eq:symmetry}),
$AJ_\beta(N,N)B\in \fra\star\frb$.
\hfill $\square$

\section{Representations of $\GLL$\label{s:5}}

\COUNTERS

Here we prove statements formulated in Subsect. \ref{ss:i-generality}, i.e., upper estimates
of the set of unitary representations of $\GLL$. The category $\cGL(\F)$ of double
cosets is an {\it ordered category} in the sense of  book \cite{Ner-book}, the statements
of Subsect. \ref{ss:i-generality}
are a kind of 'general nonsense' related to ordered categories. 

\sm

{\bf\punct Notation.}  We will regard each space $\F_\alphA$ as a space with
a distinguished basis $e_{\alpha_-+1}, \dots, e_{\alpha_+}$.
Let us assign   circles on the integer 'line' $\Z$ to basis elements of $\F_\alphA$:
$$
 \epsfbox{diagrams.1}
$$
In this section we use several linear relations $\chi(\fra):\F_\alphA\rra F_\betA$, which are spanned by vectors of the type
$$
e_i\oplus 0, e_j\oplus e_k, 0\oplus e_m. 
$$ 
We  represent such relations as two-line diagrams:

\sm

--- having a vector $e_j\oplus e_k$, we connect the corresponding circles
in the upper row and lower row; 

\sm 

--- for a vector $e_i\oplus 0$
we draw a black circle at $i$-th position in the upper row;

\sm 

--- for $0\oplus e_m$ we draw a black circle at $m$-th position of lower row.

\sm 

To draw a pair $(\chi(\fra), \eta(\fra))$ we also add to the diagram
$\eta(\fra)$ copies of the simbol $\oslash$:
$$
\epsfbox{diagrams.2}
$$
For a passage to the adjoint morphism, we must permute rows and evaluate
a new number of $\oslash$'s with Theorem \ref{th:isomorphism}.c.

\sm

{\bf \punct Proof of Lemma \ref{l:z}.\label{ss:z}} Now we must proof that {\it the number
$z$ determining a spherical character is nonnegative}.
Consider an $\alphA=(\alpha_-,\alpha_+)$ such that $H_\alphA\ne 0$.
Denote $\alphA':=(\alpha_--1,\alpha_+)$. Consider a morphism
$\frm:\alphA\to \alphA'$ defined by the diagram
$$
\epsfbox{diagrams.3}
$$
Then $\frm^*$ corresponds to the diagram
$$
\epsfbox{diagrams.4}
$$
The product $\frm^*\ast \frm$ is
$$
\epsfbox{diagrams.5},
$$
i.e.,  $\frm^*\ast \frm$ coincides with the central
element  $\zeta_\alphA^1$ defined by (\ref{eq:zeta}).
We have 
$$
0\le \rho(\frm^*)\,\rho(\frm)=\rho(\zeta_\alphA^1)
=
z(\rho)\cdot 1.
$$
 

\sm

\sm 

{\bf \punct The structure of ordered category on $\cGL(\F)$.}
Let $\betA\prec \alphA$.
Consider the morphism 
$$\lambda^\alphA_\betA:\betA\to\alphA$$ defined by
$$
\lambda^\alphA_\betA=\ov Q_\alphA\cdot 1\cdot \ov Q_\betA.
$$
The corresponding diagram has the form 
$$
\epsfbox{diagrams.6}
$$
In the notation we write a $\succ$\,-larger object $\alphA$ is superscript, and a smaller object in the subscript. 

Denote 
$$
\mu^\alphA_\betA:=\bigl(\lambda^\alphA_\betA\bigr)^*\in \Mor(\alphA,\betA).
$$
The corresponding diagram is
$$
\epsfbox{diagrams.7}
$$

Finally, define $\theta^\alphA_\betA\in\End(\alphA)$
by
$$
\theta^\alpha_\beta:=\lambda^\alpha_\beta\star  \mu ^\alpha_\beta,
$$
it corresponds to the diagram
$$
\epsfbox{diagrams.8}
$$


 Then we have
\begin{align*}
\lambda^\alphA_\betA \star \mu^\alphA_\betA&=\theta^\alphA_\betA,
\qquad
\mu^\alphA_\betA \star \lambda^\alphA_\betA=1_\betA,
\qquad
 \bigl(\theta^\alphA_\betA)^{2}=\theta^\alphA_\betA;
 \\
 (\lambda^\alphA_\betA)^*&=\mu^\alphA_\betA,\qquad
\bigl(\theta^\alphA_\betA\bigr)^*=\theta^\alphA_\betA.
\end{align*}
For
$\gamma\prec\beta\prec\alpha$ we have
\begin{equation*}
\lambda^\alphA_\betA \star \lambda^\betA_\gammA=\lambda^\alphA_\gammA,\qquad
\mu^\betA_\gammA\star \mu^\alphA_\betA  =1_\betA, \qquad
\theta^\alphA_\betA \star \theta^\alphA_\gammA=\theta^\alphA_\gammA
\end{equation*} 
This means that $\cGL(\F)$ is an ordered category with involution in the sense of
\cite{Ner-book},  Sect. III.4%
\footnote{In that definition a set of objects
	is linear ordered, but a partial order with
	existence of maximum for any pair of elements is sufficient. In any case, for Lemma \ref{l:ordered} below
	it is sufficient to consider a subcategory with two
	objects, $\betA$, $\alphA$
}. This (see \cite{Ner-book}, Lemma III.4.5, Proposition III.4.6)
implies the following statement.

\begin{lemma}
	\label{l:ordered}
	{\rm a)} Let $\betA\prec\alphA$. Then the map
	$$\iota:\frp\mapsto \mu^\alphA_\betA\star \frp\star \mu^\alphA_\betA$$
	is an embedding of semigroups $\End(\betA)\to\End(\alphA)$.
	
	\sm 
	
{\rm b)}	Let $\wh \rho$ be a $*$-representation of the category $\cGL(\F)$. 
For each object $\alphA$
denote by $H(\alphA)$ corresponding  Hilbert space.
Then the operator $\wh \rho(\lambda^\alphA_\betA):H_\betA\to H_\alpha$
is an operator of isometric embedding intertwining the representation of $\End(\betA)$
in $H(\betA)$ with the representation of $\iota\bigl(\End(\betA)\bigr)$ in the image
of the projector $\wh \rho(\theta^\alphA_\betA)$.  
\end{lemma}

{\sc Proof of Lemma \ref{l:zero-nonivertible}.b.}
Let $\deltA=(\delta_-,\delta_+)$ be a minimal element of $\Xi(\rho)$.
We must show that {\it if $\fra \in \End(\deltA)$ satisfy $\wh\rho(\fra)\ne 0$,
then $\chi(\fra)$ is an invertible matrix.}

 By Lemma \ref{l:ordered}
 for all $\epsiloN\prec \deltA$ we have $\wh\rho(\theta^\deltA_\epsiloN\bigr)=0$.
 Assume that $\rho(\fra)\ne 0$. Without loss of a generality we can assume
 $\fra^*=\fra$, otherwise we can pass to $\fra^*\star \fra$. For a  self-adjoint
 $\fra$ the linear relation $\chi(\fra)$ satisfies $\chi(\fra)^\square=\chi(\fra)$.
 Applying an appropriate conjugation by an element of $\GL$, we can reduce such a linear relation
 to a form of the type
 $$
 
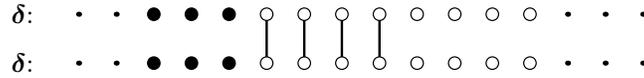
\epsfbox{diagrams.11}
 $$
 The invariant $\eta(\fra )$ can be nonzero, so $\fra$ itself has the form
 $$
 \theta^\deltA_\mU\, \zeta_\deltA^m,
 $$
 where $\mu\prec \deltA$. But $\wh \rho(\theta^\deltA_\mU )=0$
 by Lemma \ref{l:ordered}. 
   \hfill $\square$

\sm 

{\bf\punct  Proof of Theorem \ref{th:induced}.}
Recall that we wish {\it to describe possible sets 
$\Xi(\rho)$ of $\alphA$, for which the space of $\ov Q_\alphA$-fixed vectors
is non-zero.}

Let $\deltA$ be a minimal element of the set
$\Xi(\rho)$. Let $\kappA:=(\delta_-+m,\delta_++m)$, so $|\kappA|=|\deltA|$.
To be definite, assume $m>0$.
Consider the following morphism $\frr:\deltA\to\kappA$.
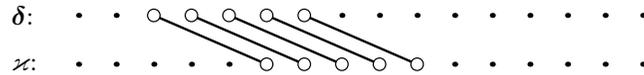
$$
\epsfbox{diagrams.9}
$$
The adjoint morphism has the form 
$$
\epsfbox{diagrams.10}
$$
(the number of $\oslash$'es is $m$).
We have
$$
\frr^*\star \frr=\zeta^m_\deltA,\qquad \frr\star \frr^*=\zeta^m_\kappA,
$$
therefore 
$$
\wh\rho(\frr)^*\wh\rho(\frr)=z^m\cdot 1,\qquad \wh\rho(\frr)\wh\rho( \frr^*)=z^m\cdot 1,
$$
If $z>0$, then the operator $z^{-m/2}\wh \rho(\frr)$ is a unitary operator
$H_\deltA\to H_\kappA$. This proves the first statement of the theorem.

\sm 

Now let $z=0$. Then $\rho(\frr)=0$. Consider a morphism $\frc:\deltA\to \kappA$.
If $\chi(\frc)$ is not a graph of an invertible operator, then 
$\chi(\frc^*\star \frc)$ also is not a graph of an invertible operator.
Therefore $\wh\rho(\frc^*\star \frc)=0$ and $\wh\rho(\frc)=0$.
If $\chi(\frc)$ is invertible, then it differs from $\chi(\frr)$ by an element
of $\GL(|\deltA|,\F)$, and $\wh\rho(\frc)=0$. Thus all operators $\rho(\frc)$
are zero, and therefore $H_\kappA=0$.

Thus (for $z=0$), if $\deltA$  a  minimal element of the set
$\Xi(\rho)$, then for $\nU\in \Xi(\rho)$ we have
$\nU=\deltA$ or $|\nU|>|\deltA|$. In particular, $\Xi(\rho)$
contains a unique minimal element. 
This proves the second statement of Theorem.

\sm

{\bf \punct The group $\GLL$ has type $I$.%
\label{ss:proof-I}}

\begin{lemma} 
	\label{l:linear-ordering}
There exists a {\rm(}noncanonical{\rm )} linear order $\vartriangleleft$ on $\cA$ compatible
with the partial order $\prec$
satisfying the condition:
 for each $\alphA\in\cA$  the set of all $\betA\in \cA$ such that  
$\betA\vartriangleleft\alphA$, is finite. 
\end{lemma}

\begin{figure}
	$$\mathrm a)\quad\epsfbox{order.1}\qquad \mathrm b)\quad\epsfbox{order.2}$$
	
	\caption{Ref. to proof of Lemma \ref{l:linear-ordering}.}
\end{figure}

{\sc Proof of Lemma.}
For $\alphA\in \cA$ denote by 
$m:=\alpha_-$, $n:=\alpha_+-\alpha_-$. Then $(m,n)$ ranges in the set
$\Z\times \Z_+$. A set $\betA\prec \alphA$ is drawn
on Fig. 1.a. Now define a linear order $\vartriangleleft$
in the following way. We consider the sequence of segments $I_0$, $I_1$, $I_2$, \dots
as it is drawn on Fig. 1.b and enumerate integer  points of the upper half-plane
in the following way: the unique point of $I_0$, then we pass $I_1$ in upper direction, 
$I_2$ in upper direction, etc.
\hfill $\square$

\sm 

{\sc Proof of Theorem \ref{th:I}.}
We must examine the von Neumann algebra $\frN$ of all operators
 commuting with all operators $\rho(g)$, where
$g\in \GLL$.

Keeping in mind Lemma \ref{l:linear-ordering},
we write a sequence
$$
\alphA_1\vartriangleleft \alphA_2 \vartriangleleft \alphA_3 \vartriangleleft\dots
$$
containing all $\alphA\in \cA$. Let us decompose the Hilbert space $H$
into a countable direct sum $H=\oplus_{j=1}^\infty K_j$ according the following inductive rule.
Consider the subspace $H_{\alphA_1}$ and its $\GLL$-cyclic span $K_1$. Next, consider
$H\ominus K_1$, the subspace $(H\ominus K_1)_{\alphA_2}$ and its $\GLL$-cyclic span $K_2$.
Then we consider the cyclic span  $K_2$ of $(H\ominus K_1\ominus K_2)_{\alphA_3}$. Etc.

Clearly, elements of $\frN$ leave all subspaces $K_j$ invariant,
and $\frN=\oplus \frN_j$, where $\frN_j$ are induced von Neumann algebras
in $K_j$. Therefore, we must examine $\frN_j$. It is easy to see that
this algebra is isomorphic to the algebra of operators in $(K_j)_{\alphA_j}$
commuting with $\End(\alphA_j)$. The latter semigroup is  
 $\GL(|\alphA_j|,\F)\times \Z_+$.
 Evidently the von Neumann algebra generated by this semigroup
 has type I, therefore its commutant $\frN_j$ also has type I.
 This proves the statement a) of the theorem (the group $\GLL$ has type I).
 
 The statement b) follows from the same considerations. It is clear
 that any $*$-presentation of $\GL(|\alphA_j|,\F)\times \Z_+$
 in $H_{\alphA_j}$
 can be decomposed into a direct integral. A simple watching shows that this
 induces a decomposition of the whole space $K_j$ into a direct integral. 
\hfill $\square$

\sm

{\bf \punct Constructions of all representations of $\GLL$ with $z=0$.%
	\label{ss:constructions}} 
A proof of Proposition \ref{pr:z=0} is an exercise on induced representations.
The homogeneous space 
$$X:=\ov P_\alphA\setminus \GLL\simeq P_\alphA\setminus \GL(\infty,\F)$$
is discrete
(it  consists of two-terms flags of the form $Y\supset Z$, where
$Y\subset \V$ is a compact subspace of volume $q^{\alpha_-}$ and
$Z$ is a compact subspace of volume  $q^{\alpha_+}$).
This allows to apply the usual construction of induced representations, see \cite{Mack-induced},
\cite{Kir}, Subsect. 13.1.

\sm

Namely, let $\tau$ be an irreducible representation of the finite group
$\GL(|\alphA|,\F)$ in the space $L_\tau$, we can consider $\tau$ as a representation of
$\ov P_\alphA$ trivial on $\ov Q_\alphA$.
Consider the representation $T$ of $\GLL$ induced from the representation $\tau$ 
of $\ov P_\alphA$ (see, e.g., \cite{Kir}, Subsect. 13.1). It is realized in the space $\ell^2(X,L_\tau)$
of $L_\tau$-valued  functions on $X$ with the inner product
$$
\la F_1, F_2\ra_{\ell^2(X,L_\tau)}:=
\sum_{x\in X}\la F_1(x),F_2(X)\ra_{L_\tau}.
$$
The group $\GLL$ acts in $\ell^2(X,L_\tau)$ by transformations
of the form
$$
T(g)F(x):=U(g,x)F(xg),
$$
where $U(g,x)$ is a function on $\GLL\times X$ taking values in
the unitary group of the space $L_\tau$ and satisfying the conditions
$$
U(g_1g_2,x)= U(x,g_1)\,U(xg_1,g_2); \qquad U(1_{\GLL},x)=1_{L_\tau}
$$ 
(under this conditions $T(g)$ is a representation) and
$$
U(h,x_0)=\tau (h),\qquad\text{where $x_0\in X$ is the initial point
	and $h\in \ov P_\alphA$}.
$$
(this determines a function $U(g,x)$ up to a natural equivalence).
The usual  Mackey arguments \cite{Mack-induced} show that the
representation $T(g)$ is irreducible. To a verification of irreducibility
we need only two facts:

\sm 

1) $V_\tau$ are finite-dimensional;

\sm 

2) all orbits of $Q_\alphA$ on $Q_\alphA\setminus \GL(\infty,\F)$
except the point $x_0$ are infinite.

\sm 

Formally, we can refer to \cite{Cor}, Theorem 2.

The property 2) easily implies also that the subspace of $Q_\alphA$-fixed vectors 
in $\ell^2(X,L_\tau)$ 
consists of functions supported by one point $x_0$. I.e., this space is $L_\tau$.





 \tt
\noindent
Yury Neretin\\
Wolfgang  Pauli Institute/c.o. Math. Dept., University of Vienna \\
\&Institute for Theoretical and Experimental Physics (Moscow); \\
\&MechMath Dept., Moscow State University;\\
\&Institute for Information Transmission Problems;\\
yurii.neretin@univie.ac.at
\\
URL: http://mat.univie.ac.at/$\sim$neretin/

\end{document}